\documentclass[a4paper]{article}
\usepackage{a4wide}

\usepackage {tikz}
\usetikzlibrary {positioning}
\usepackage[english]{babel}
\usepackage[utf8]{inputenc}
\usepackage{ifthen}
\usepackage{color}
\usepackage{amsmath,amsthm,amssymb,amsfonts}
\usepackage{ stmaryrd }

\usepackage{bbold}
\usepackage{ mathrsfs }
\usepackage{bbm}
\usepackage{graphicx}
\usepackage{psfrag}
\usepackage[normalem]{ulem}
\usepackage{enumitem}   
\usepackage{relsize}
\usepackage{subfig}

\usepackage[makeroom]{cancel}

\usepackage[unicode,colorlinks=true,pagebackref=false]{hyperref}

\definecolor{change}{rgb}{0,.55,.55}

\newtheorem{Def}{Definition}
\newtheorem{Thm}[Def]{Theorem}
\newtheorem{Rmk}[Def]{Remark}
\newtheorem{Lemma}[Def]{Lemma}
\newtheorem{Corollary}[Def]{Corollary}

\newtheorem{Not}[Def]{Notation}
\newtheorem{Ass}[Def]{Assumption}

\newcommand{\R}{\mathbb{R}}

\newcommand{\N}{\mathbb{N}}

\newcommand{\Cspace}{C}

\newcommand{\generator}{L}

\newcommand{\Kalman}{K}
\newcommand{\aShift}{a}

\newcommand\indicator{\mathbb{1}}

\newcommand\norm[2]{\left\|#1\right\|_{#2}}
\newcommand\comCov[2]{\left\llbracket #1 | #2\right\rrbracket}

\title{Analysis of the Ensemble Kalman--Bucy Filter for correlated observation noise}
\author{Sebastian Ertel, Wilhelm Stannat}
\date{}

\begin{document}

	\maketitle	
	
	\textbf{Abstract.} Ensemble Kalman--Bucy filters (EnKBFs) are an important tool in Data Assimilation that aim to approximate the posterior distribution for continuous time filtering problems using an ensemble of interacting particles. In this work we extend a previously derived unifying framework for consistent representations of the posterior distribution to correlated observation noise and use these representations to derive an EnKBF suitable for this setting as a constant gain approximation of these optimal filters. Existence and uniqueness results for both the EnKBF and its mean field limit are provided. The existence and uniqueness of solutions to its limiting McKean-Vlasov equation does not seem to be covered by the existing literature.
	In the correlated noise case the evolution of the ensemble depends also on the pseudoinverse of its empirical covariance matrix, which has to be controlled for global well posedness. These bounds may also be of independent interest.
	Finally the convergence to the mean field limit is proven. The results can also be extended to other versions of EnKBFs.
	
	\tableofcontents
	
	\section{Introduction}
	
	The aim of filtering algorithms is to approximate the state of a signal from noisy and potentially incomplete or indirect observations. Ensemble Kalman--Bucy filters\footnote{The time continuous version of classical Ensemble Kalman filters.} (EnKBFs)  accomplish this by employing a system of interacting particles, which then are propagated according to the signal dynamics enriched by an added interaction term that nudges them in the direction of the observations.\newline
	
	Ensemble Kalman(--Bucy) Filters were first introduced in \cite{Evensen} and are nowadays widely used for data assimilation tasks in many scientific fields such as meteorology and the geosciences. For uncorrelated observation noise there exists a rich literature treating the mathematical theory of EnKBFs. See for example the recent overview papers \cite{On the mathematical theory of linear EnKBF},\cite{ReichStuart}, and the references found therein. In this work we derive an EnKBF for correlated observation noise and analyse its mean field limit.\newline
	
	More precisely for an arbitrary, but fixed timeframe $T>0$ we consider a $\R^{d_x}$-valued signal process $\left(X_t\right)_{t\in[0,T]}$ determined by the stochastic differential equation (SDE)
	\begin{align}\label{SDE correlated observation noise}
		\mathrm{d}X_t=
		B_t(X_t)\mathrm{d}t+C_t(X_t)\mathrm{d}W_t+\tilde{C}_t(X_t)\mathrm{d}V_t.
	\end{align}
	
	Both $W$ and $V$ shall be two independent $\R^{d_w}$- and $\R^{d_v}$-dimensional Brownian motions. Furthermore we make the following standard assumption regarding the coefficients of the equation.
	
	\begin{Ass}\label{standard assumptions}
		The drift $B:[0,T]\times\R^{d_x}\to\R^{d_x}$ and the (square root of the) diffusion $C:[0,T]\times\R^{d_x}\to\R^{d_x\times d_w}$ are Borel-measurable and satisfy the usual  linear growth and global Lipschitz conditions found for example in \cite{KaratzasShreve}. The corresponding Lipschitz constants shall be denoted by $\mathrm{Lip}(B)$ and $\mathrm{Lip}(C)$. $\tilde{C}:[0,T]\times\R^{d_x}\to\R^{d_x\times d_v}$ is assumed to be a continuous matrix valued function. 
	\end{Ass}

	The $R^{d_y}$-valued observation process $\left(Y_t\right)_{t\in[0,T]}$ is then given by
	\begin{align*}
		\mathrm{d}Y_t=H_t X_t\mathrm{d}t+\Gamma_t\mathrm{d}V_t,
	\end{align*}
	where both $H:[0,T]\to\R^{d_y\times d_x}$ and $\Gamma:[0,T]\to\R^{d_y\times d_v}$ are continuous. As is usual, we assume that $Y_0=0$.\newline
	
	The goal of stochastic filtering is to compute, or at least approximate, the conditional distribution of $X_t$ given all past observations $(Y_s)_{s\in[0,t]}$, which we denote by
	\begin{align}\label{definition posterior}
		\eta_t:=\mathbb{P}\left(~X_t\in\cdot~|~Y_{s},~s\leq t~\right),
	\end{align}
	for all times $t\in[0,T]$. $\eta$ is referred to as the posterior distribution or the optimal filter. We shall denote the integral of a testfunction $\phi$ with respect to $\eta_t$ by $\eta_t(\phi):=\int_{\R^{d_x}} \phi(x)~\eta_t(\mathrm{d}x)$.\newline
	
	For any sufficiently regular and integrable testfunction $\phi$, the weak form of the Kushner--Stratonovich equation \cite{BainCrisan}, \cite{NueReiRoz}
	\begin{align}\label{correlated KSE}
		\begin{split}
			\mathrm{d}\eta_t(\phi)
			&=
			\eta_t(\generator_t\phi)\mathrm{d}t
			\\&\phantom{=}+
			\left(\eta_t\left(H\mathrm{id}_{\R^{d_x}}\phi\right)^{\mathrm{T}}
			-\eta_t\left(H_t\mathrm{id}_{\R^{d_x}}\right)^{\mathrm{T}}\eta_t\left(\phi\right)
			+\eta_t\left((\nabla\phi)^{\mathrm{T}}\tilde{C}_t\Gamma^{\mathrm{T}}_t\right)
			\right)
			R^{-1}_t \mathrm{d} I_t
			\\
			\mathrm{d} I_t&=
			\mathrm{d}Y_t-\eta_t(H_t \mathrm{id}_{\R^{d_x}})\mathrm{d}t
		\end{split}
	\end{align}
	is satisfied. Hereby $\generator_t$ is the generator associated to \eqref{SDE correlated observation noise} and  $R_t:=\Gamma_t\Gamma_t^{\mathrm{T}}$. In many applications the dimension of the signal $d_x$ is too large, so that an approximation of \eqref{correlated KSE} using standard numerical PDE solvers becomes infeasible and even in one dimension they are not efficient for unstable/transient signals. Sequential Monte-Carlo type methods called particle filters provide an alternative approximation of the optimal filter that is often considered, however these methods can suffer from weight degeneracy which makes them scale poorly with the state dimension.  Instead practitioners often rely on EnKBFs, which, even though they only provide an approximation to the optimal filter in the linear Gaussian case, have proven to be a successful tool for Data Assimilation tasks even for nonlinear signals. A particular version of EnKBF that will be the main focus of this paper is
	\begin{align}\label{nonlinear Ensemble Kalman - introduction}
		\begin{split}
			\mathrm{d}X^i_t
			&=
			B_t\left(X^i_t\right)\mathrm{d}t+
			C_t\left(X^i_t\right)\mathrm{d}W^i_t
			+
			\tilde{C}_t\mathrm{d}V^i_t
			\\&\phantom{=}+
			\left(P^M_t H^\mathrm{T}_t+\tilde{C}_t\Gamma^{\mathrm{T}}_t\right)
			R^{-1}_t
			\left(\mathrm{d}Y_t-\frac{H_t\left(X^i_t+x^M_t\right)}{2}\mathrm{d}t\right)
			\\&\phantom{=}-
			\left(	P^M_t H^\mathrm{T}_t+\tilde{C}_t\Gamma^{\mathrm{T}}_t\right)~
			R^{-1}_t~
			\Gamma_t\tilde{C}^{\mathrm{T}}_t\left(P^M_t\right)^{+}
			\frac{X^i_t-x^M_t}{2}\mathrm{d}t
		\end{split}
	\end{align}
	for $i=1,\cdots,M$, with
	\begin{align*}
		x^M_t:=\frac{1}{M}\sum_{i=1}^{M} X^i_t
		~\text{and}~
		P^M_t:=\frac{1}{M-1}\sum_{i=1}^{M}\left(X^i_t-x^M_t\right)\left(X^i_t-x^M_t\right)^\mathrm{T},
	\end{align*}
	and where $\left(P^M_t\right)^{+}$ denotes the Moore--Penrose pseudoinverse of $P^M_t$. For $\tilde{C}=0$ equation \eqref{nonlinear Ensemble Kalman - introduction} is the continuous-time counterpart to the deterministic Ensemble Kalman filter introduced in \cite{Sakov Oke}. We will prove well posedness of both \eqref{nonlinear Ensemble Kalman - introduction} and its mean field limit and also show propagation of chaos. Our analysis can easily be extended to other types of EnKBFs.\newline

	Let us now briefly summarize the structure  and the contributions of this paper:
	\begin{itemize}
		\item In section \ref{section mf representation} we formally derive a McKean--Vlasov equation that is consistent with \eqref{correlated KSE}, meaning that if a solution exists, its law will evolve according to \eqref{correlated KSE}. For the uncorrelated case $\tilde{C}=0$ this has already been done in \cite{PathirajaStannatReich} in a unifying framework covering existing mean field filters like  \cite{CrisanXiong} and the well known feedback particle filter (FPF) \cite{YangMehtaMeyn1}, \cite{YangMehtaMeyn2}. In the correlated noise case FPFs have already been derived in \cite{LuoMiao} and \cite{NueReiRoz}. In this work we extend the framework in \cite{PathirajaStannatReich} to arbitrary $\tilde{C}$.\newline
		
		\item In section \ref{EnKBF section} we derive a mean field EnKBF as a special case of the consistent mean field representation of section \ref{section mf representation} in a linear Gaussian filtering setting. In the uncorrelated noise framework continuous time McKean--Vlasov interpretations of the filtering equation with Ensemble Kalman filtering mean field interpretations for linear Gaussian models have been discussed for the first time in \cite{DelMoralTugaut}. We use arguments first found in \cite{DelMoralTugaut} to prove the well posedness of the mean field EnKBF in the correlated noise framework.\newline
		
		\item In section \ref{EnKBF for nonlinear signals - section} we first interpret the mean field EnKBF in the nonlinear, non Gaussian setting as an approximation of the optimal filter derived in section \ref{section mf representation}. In the uncorrelated setting this connection is often referred to as the constant gain approximation of the feedback particle filter \cite{TaghvaeiDeWiljesEtAl}.\\
		In subsection \ref{section proof of thm} we then show the well posedness of the mean field EnKBF for nonlinear signals and correlated observation noise. Previous results on this matter \cite{Coghi et al} require the observation function $H$ to be bounded. Linear observations, which are highly relevant in practice, do not seem to be covered by existing literature, even in the simpler case of uncorrelated observation noise. We prove the well posedness for this case by using a mixture of a fixed point and a deterministic localization argument.\\
		Finally, subsection \ref{section other EnKBF} discusses other mean field EnKBFs. Our well posedness result can be extended to these McKean--Vlasov equations in a straightforward manner.\newline
		
		\item In section \ref{section particle system} we show the well posedness of of the particle system approximating the mean field EnKBF, which we just refer to as the EnKBF. For uncorrelated observation noise this was already proven in \cite{StannatLange}, however the  correlated noise framework requires an additional term depending on the (pseudo) inverse of the ensemble covariance matrix, that may become singular. We prove existence and uniqueness of strong solutions under suitable assumptions by showing that these singularities will never be hit by the filter. Due to the stochasticity of the ensemble covariance matrix we cannot apply techniques for bounding the spectrum from below found for example in \cite{DeWiljesReichStannat} and instead dominate its inverse in terms of  the ensemble covariance itself. This technique may be of independent interest.\newline
		
		\item In section \ref{section propagation of chaos} we derive a time dependent propagation of chaos result, extending previous results in the uncorrelated case found for example in \cite{StannatLange}. In the uncorrelated, linear Gaussian setting propagation of chaos has first been proven in \cite{DelMoralTugaut}. In this setting the mean field limit is the optimal filter and several uniform in time estimates for the convergence exist (e.g.\cite{Riccati diff - stability},\cite{Riccati diff - Perturbation},\cite{DelMoralTugaut},\cite{DelMoralKurtzmannTugaut}), providing error/inconsistency bounds for the EnKBF. In the nonlinear setting our propagation of chaos result no longer gives these error bounds as the mean field limit does not coincide with the posterior. Nevertheless the mean field limit provides a connection to the optimal filter via the constant gain approximation discussed in subsection \ref{subsection constant gain approximation} and as such may be an important tool for quantifying the inconsistency of the EnKBF.\newline
		
	\end{itemize}

	\section{Mean field representation of the posterior}\label{section mf representation}	
	
	\begin{Ass}
		To achieve greater generality in our results throughout this section, we allow both $\tilde{C}$ and $H$ to be state dependent functions. I.e. for every $t\geq 0$ we consider the (nonlinear) maps $x\mapsto\tilde{C}_t(x)$ and $x\mapsto H_t(x)$. This assumption is restricted to section \ref{section mf representation}.
	\end{Ass}
	
	{In \cite{PathirajaStannatReich} a representation of the posterior distribution in the uncorrelated setting by a diffusion process was derived, by matching the Kolmogorov equation of said diffusion to the Kushner--Stratonovich equation. This gives a class of mean field equations that contains many well known optimal filters like the Feedback Particle Filter \cite{YangMehtaMeyn1}. In this section we follow the approach of \cite{PathirajaStannatReich} in order to derive consistent mean--field representations in the correlated noise framework. Thus we aim to represent the posterior through a McKean--Vlasov equation}, i.e. we want to find a diffusion process $\left(\bar{X}_t\right)_{t\geq 0}$ such that for all times $t>0$, its (conditional) marginal law $\bar{\eta}_t$ is given by the posterior $\eta_t$. To this end we assume throughout this section that all functions appearing are sufficiently regular and integrable, so that we can always differentiate and integrate whenever necessary.\newline

	Let $\bar{W}$ be an independent copy of the Brownian motion $W$ and $\bar{V}$ be an independent copy of $V$. To determine $\bar{X}$ we make the Ansatz
	\begin{align}\label{mean field ansatz}
		\begin{split}
			\mathrm{d}\bar{X}_t
			&=B_t(\bar{X}_t)\mathrm{d}t
			+C_t(\bar{X}_t)\mathrm{d}\bar{W}_t
			+\tilde{C}_t(\bar{X}_t)\mathrm{d}\bar{V}_t
			\\&\phantom{=}+\Kalman_t\left(\bar{X}_t,Y_{0:t},\bar{\eta}_t\right)\mathrm{d}Y_t
			+\aShift_t\left(\bar{X}_t,Y_{0:t},\bar{\eta}_t\right)\mathrm{d}t,
		\end{split}
	\end{align}
	for two functions $\Kalman$ and $\aShift$ depending on both the state and the (marginal) law of the process $\bar{X}$. To make the class of potential processes $\bar{X}$ even larger, we also allow for a dependence on all past observations $Y_{0:t}$, which as we will later on see, can be dropped. Furthermore let us make the following notational convention, which helps to shorten formulas.
	
	\begin{Not}\label{definition EY}
		For any random variable $Z$ we denote the conditional expectation with respect to all observations $Y$ on the total time interval $[0,T]$ by
		\begin{align*}
			\mathbb{E}_{Y}\left[Z\right]
			:=
			\mathbb{E}\left[~Z~|~Y_t,~t\leq T~\right].
		\end{align*}
	\end{Not}

	\begin{Rmk}
		Notation \ref{definition EY} is in accordance to \cite{PathirajaStannatReich}. Due to the independence of $Y$ from $\bar{W}$ and $\bar{V}$, another suitable way to define $\mathbb{E}_{Y}$ would be to let it denote the integral with respect to the joint law by $\bar{W}$ and $\bar{V}$. In both cases we have
		\begin{align*}
			\mathbb{E}\left[Z\right]=\int \mathbb{E}_{Y}\left[Z\right]~\mathbb{P}(Y\in\mathrm{d}y).
		\end{align*}
	\end{Rmk}
	
	\begin{Rmk}
		{Another way to  view the McKean--Vlasov equations derived in this and the subsequent sections is to interpret the observations $Y$ as a deterministic rough path. This is a research direction that has recently received attention in the EnKBF literature (e.g. \cite{Coghi et al}) and is a natural modelling choice, as in practice one will often only have the data of a single realization/path of $Y$.} Large parts of our calculations, in particular all results in the sections \ref{section mf representation}, \ref{EnKBF section} and \ref{EnKBF for nonlinear signals - section}, also hold in this setting. In particular the well posedness proof in subsection \ref{section proof of thm} can easily be addapted to the case when $Y$ is a deterministic rough path, using the recent well posedness result for rough-stochastic differential equations found in \cite{Friz}.\newline
	\end{Rmk}
	
	Let $\phi:\R^{d_x}\to\R$ be some sufficiently regular and bounded testfunction. Using Itô's formula and taking the (conditional) expectation $\mathbb{E}_Y$, we derive for $\bar{\eta}_t(\phi):=\mathbb{E}_{Y}\left[\phi\left(\bar{X}_t\right)\right]$ the equation
	\begin{align}\label{Fokker Planck for MF Filter}
		\begin{split}
			\mathrm{d}\bar{\eta}_t(\phi)
			&=\bar{\eta}_t(\generator_t\phi)\mathrm{d}t+
			\bar{\eta}_t\left(\nabla\phi\cdot\Kalman_t\left(\cdot,Y_{0:t},\bar{\eta}_t\right)\right)\mathrm{d}Y_t
			\\&\phantom{=}+\bar{\eta}_t\left(\nabla\phi\cdot \aShift_t\left(\cdot,Y_{0:t},\bar{\eta}_t\right)\right)\mathrm{d}t
			\\&\phantom{=}+
			\frac{1}{2}\bar{\eta}_t\left(\mathrm{tr}\left[\phi'' \Kalman_t\left(\cdot,Y_{0:t},\bar{\eta}_t\right) R_t \Kalman_t\left(\cdot,Y_{0:t},\bar{\eta}_t\right)^{\mathrm{T}}\right]\right)
			\mathrm{d}t.
		\end{split}
	\end{align}
	Hereby $\nabla\phi$ denotes the gradient of $\phi$ and $\phi''$ its Hessian.\newline
	
	Since $\bar{\eta}_t$ shall coincide with the posterior $\eta_t$, it must also adhere to the Kushner--Stratonovich equation \eqref{correlated KSE}. By comparing the terms on the right-hand side of both equations \eqref{correlated KSE} and \eqref{Fokker Planck for MF Filter}, we get the two consistency conditions
	\begin{align}\label{consistency K}
		\begin{split}
			&\bar{\eta}_t\left(\nabla\phi\cdot\Kalman_t\left(\cdot,Y_{0:t},\bar{\eta}_t\right)\right)
			=
			\left(\bar{\eta}_t\left(H_t^{\mathrm{T}}\phi\right)
			-\bar{\eta}_t\left(H_t^{\mathrm{T}}\right)\bar{\eta}_t\left(\phi\right)
			+\bar{\eta}_t\left((\nabla\phi)^{\mathrm{T}}\tilde{C}_t\Gamma^{\mathrm{T}}_t\right)
			\right) R^{-1}_t
		\end{split}
	\end{align}
	and
	\begin{align}\label{consistency a}
		\begin{split}
			&\bar{\eta}_t
			\left(\nabla\phi\cdot \aShift_t\left(\cdot,Y_{0:t},\bar{\eta}_t\right)\right)
			+
			\frac{1}{2}\bar{\eta}_t\left(\mathrm{tr}\left[\phi'' \Kalman_t\left(\cdot,Y_{0:t},\bar{\eta}_t\right) R_t
			\Kalman_t\left(\cdot,Y_{0:t},\bar{\eta}_t\right)^{\mathrm{T}}\right]\right)
			\\&=
			-\left(\bar{\eta}_t\left(H_t^{\mathrm{T}}\phi\right)
			-\bar{\eta}_t\left(H_t^{\mathrm{T}}\right)\bar{\eta}_t\left(\phi\right)
			+\bar{\eta}_t\left((\nabla\phi)^{\mathrm{T}}\tilde{C}_t\Gamma^{\mathrm{T}}_t\right)
			\right) R^{-1}_t \bar{\eta}_t\left(H_t\right).
		\end{split}
	\end{align}
	
	Since the right hand sides of both equations do not depend on past observations $Y_{0:t}$, we can drop the $Y$-dependence in both $\Kalman$ and $\aShift$. Note that both $ \Kalman_t\left(\cdot,\bar{\eta}_t\right)$ and $\aShift_t\left(\cdot,\bar{\eta}_t\right)$ are thus purely statistical quantities of the distribution $\bar{\eta}_t$.\newline
	
	Assuming that $\bar{\eta}_t$ admits a sufficiently regular density, which shall also be denoted by $\bar{\eta}_t$, we derive more direct characterizations of the two terms $\Kalman$ and $\aShift$ in the following.\newline

	\subsection{Consistency of the Kalman-gain}
	
	First we investigate the consistency condition \eqref{consistency K} for the Kalman gain term $\Kalman$. To this end we make the following notational conventions that will be used throughout this paper. 
	
	\begin{Not}
		The divergence of a matrix-valued function $A$ shall be interpreted columnwise, i.e. if $A(x)\in\R^{m\times k}$ we have
		\begin{align*}
			\mathrm{div}(A)
			:=\left(\mathrm{div}(A_{\cdot,1}),\cdots,\mathrm{div}(A_{\cdot,k})\right)
			:=\left(\sum_{j=1}^{m}\partial_{x_j}A_{j,1},\cdots,\sum_{j=1}^{m}\partial_{x_j}A_{j,k}\right).
		\end{align*}
		Furthermore we also interpret the scalar product between a matrix and a vector columnwise. Thus the scalar product between $A$ and the gradient $\nabla$ gives the following row vector-valued differential operator
		\begin{align*}
			A\cdot\nabla
			:=\left(A_{\cdot,1}\cdot\nabla,\cdots,A_{\cdot,k}\cdot\nabla\right)
			:=\left(\sum_{j=1}^{m}A_{j,1}\partial_{x_j},\cdots,\sum_{j=1}^{m} A_{j,k}\partial_{x_j}\right).
		\end{align*}
		
	\end{Not}
	
	Employing integration by parts and the Fundamental Theorem of Calculus of Variations, we see that \eqref{consistency K} can be interpreted as the weak form of the partial differential equation
	\begin{align}\label{consistency K - differential form}
		-\mathrm{div}\left(\bar{\eta}_t\Kalman_t(\cdot,\bar{\eta}_t)\right)+\mathrm{div}\left(\bar{\eta}\tilde{C}_t\Gamma^{\mathrm{T}}_t\right)R^{-1}_t
		=
		\left(H^{\mathrm{T}}-\bar{\eta}_t(H^{\mathrm{T}})\right)R^{-1}_t\bar{\eta}_t.
	\end{align}
	
	If we denote the dependence of $\Kalman$ on $\tilde{C}$ by $\Kalman^{\tilde{C}}$, then we clearly have
	\begin{align}\label{K as a translation}
		\Kalman^{\tilde{C}}_t=\Kalman^{0}_t+\tilde{C}_t\Gamma^{\mathrm{T}}_t R^{-1}_t,
	\end{align}
	with
	\begin{align}\label{K0}
		-\mathrm{div}\left(\bar{\eta}_t\Kalman^{0}_t(\cdot,\bar{\eta}_t)\right)
		=
		\left(H_t^{\mathrm{T}}-\bar{\eta}_t(H_t^{\mathrm{T}})\right)R^{-1}_t\bar{\eta}_t.
	\end{align}
	
	Thus $\Kalman^{\tilde{C}}$ is just a translation of $\Kalman^{0}$, the Kalman gain for the uncorrelated case, which in turn is defined uniquely up to $\bar{\eta}_t$-harmonic vector fields.\newline

	\begin{Rmk}
		By using partial integration to rewrite \eqref{K0} into flux form
		\begin{align*}
			\int_{\partial D} \bar{\eta}_t\Kalman^{0}_t(\cdot,\bar{\eta}_t)\cdot(-\nu_D)\mathrm{d}s
			=
			\int_{D} \left(H_t^{\mathrm{T}}-\bar{\eta}_t(H_t^{\mathrm{T}})\right)R^{-1}_t\bar{\eta}_t\mathrm{d}x,
		\end{align*}
		where $D$ is an arbitrary domain and $\nu_D$ is its outer normal vector, we see that $\Kalman^{0}_t(\cdot,\bar{\eta}_t)$ can be interpreted as the velocity (speed and direction) by which particles of density $\bar{\eta}_t$ must travel, such that for every domain $D$ the flux into $D$ is equal to the difference between the expected observation in that domain and the global expected observation. Therefore particles are pushed into areas where the observations are expected to deviate from average observation in the whole ensemble. Thus $\Kalman^{0}_t(\cdot,\bar{\eta}_t)$ can be seen as a quantity of the distribution $\bar{\eta}_t$ that corresponds to the exploration of the state space with respect to the observation function $H$.
	\end{Rmk}

	\subsection{Consistency of the correctional transport term}
	
	Next we investigate the consistency equation \eqref{consistency a} for the correctional transport term $a$. Note that while $\Kalman$ also shows up in this equation, it is already fully determined by \eqref{consistency K}.\newline

	By again using integration by parts and the Fundamental Theorem of Calculus of Variations, we derive the strong form of consistency condition \eqref{consistency K}
	\begin{align}\label{consistency a - differential form 1}
		\begin{split}
			&-\mathrm{div}\left(\bar{\eta}_t\aShift_t\left(\cdot,\bar{\eta}_t\right)\right)
			+
			\frac{1}{2}
			\sum_{i,j=1}^{d_x}
			\partial^2_{x_i x_j}
			\left(\bar{\eta}_t K_t\left(\cdot,\bar{\eta}_t\right) R_t K_t\left(\cdot,\bar{\eta}_t\right)^{\mathrm{T}}\right)_{ji}
			\\&\phantom{=}-
			\mathrm{div}\left(\bar{\eta}_t\tilde{C}_t\Gamma^{\mathrm{T}}_t\right)
			R^{-1}_t \bar{\eta}_t\left(H_t\right)
			=
			-\left(H^{\mathrm{T}}_t-\bar{\eta}_t(H^{\mathrm{T}}_t)\right)R^{-1}_t \bar{\eta}_t(H_t)\bar{\eta}_t
			.
		\end{split}
	\end{align}

	%

	Using
	\begin{align*}
		&\sum_{i,j=1}^{d_x}
		\partial^2_{x_i x_j}
		\left(\bar{\eta}_t \Kalman_t\left(\cdot,\bar{\eta}_t\right) R_t 
		\Kalman_t\left(\cdot,\bar{\eta}_t\right)^{\mathrm{T}}\right)_{ji}
		\\&=
		\mathrm{div}\left(
		\Kalman_t\left(\cdot,\bar{\eta}_t\right)~R_t~
		\mathrm{div}\left(\bar{\eta}_t \Kalman_t\left(\cdot,\bar{\eta}_t\right)\right)^{\mathrm{T}}
		\right)
		+
		\mathrm{div}\left(\bar{\eta}_t \left(
		\left(\Kalman_t\left(\cdot,\bar{\eta}_t\right)\cdot\nabla\right)
		R_t \Kalman_t^{\mathrm{T}}\left(\cdot,\bar{\eta}_t\right)
		\right)^{\mathrm{T}}\right)
	\end{align*}
	and the consistency equation of the Kalman gain term \eqref{consistency K - differential form} gives us the identity
	\begin{align*}
		&\sum_{i,j=1}^{d_x}
		\partial^2_{x_i x_j}
		\left(\bar{\eta}_t \Kalman_t\left(\cdot,\bar{\eta}_t\right) R_t 
		\Kalman_t\left(\cdot,\bar{\eta}_t\right)^{\mathrm{T}}\right)_{ji}
		\\&=
		-\mathrm{div}\left(\bar{\eta}_t K_t\left(\cdot,\bar{\eta}_t\right)\left(H_t-\bar{\eta}_t(H_t)\right)\right)
		+
		\mathrm{div}\left( K_t\left(\cdot,\bar{\eta}_t\right)\mathrm{div}\left(\bar{\eta}_t\tilde{C}_t\Gamma^{\mathrm{T}}_t\right)^{\mathrm{T}}\right)
		\\&\phantom{=}+
		\mathrm{div}\left(\bar{\eta}_t\left(
		\left(\Kalman_t\left(\cdot,\bar{\eta}_t\right)\cdot\nabla\right)
		R_t \Kalman_t^{\mathrm{T}}\left(\cdot,\bar{\eta}_t\right)
		\right)^{\mathrm{T}}\right).
	\end{align*}
	
	Thus we can rewrite \eqref{consistency a - differential form 1} into
	\begin{align*}
		\begin{split}
			&-\mathrm{div}\left(\bar{\eta}_t\aShift_t\left(\cdot,\bar{\eta}_t\right)\right)
			-\frac{\mathrm{div}\left(\bar{\eta}_t \Kalman_t\left(\cdot,\bar{\eta}_t\right)\left(H_t-\bar{\eta}_t(H_t)\right)\right)}{2}
			-
			\mathrm{div}\left(\bar{\eta}_t\tilde{C}_t\Gamma^{\mathrm{T}}_t\right)R^{-1}_{t}\bar{\eta}_t\left(H_t\right)
			\\&\phantom{=}+
			\frac{\mathrm{div}\left( \Kalman_t\left(\cdot,\bar{\eta}_t\right)\mathrm{div}\left(\bar{\eta}_t\tilde{C}_t\Gamma^{\mathrm{T}}_t\right)^{\mathrm{T}}\right)}{2}
			+
			\frac{\mathrm{div}\left(\bar{\eta}_t
				\left(
				\left(\Kalman_t\left(\cdot,\bar{\eta}_t\right)\cdot\nabla\right)
				R_t \Kalman_t^{\mathrm{T}}\left(\cdot,\bar{\eta}_t\right)
				\right)^{\mathrm{T}}\right)}{2}
			\\&=
			-\left(H_t^{\mathrm{T}}-\bar{\eta}_t(H_t^{\mathrm{T}})\right)R^{-1}_{t}\bar{\eta}_t(H_t)\bar{\eta}_t.
		\end{split}
	\end{align*}

	Just as in	\cite{PathirajaStannatReich} one sees immediately that
	{\small\begin{align}\label{a formula}
			\begin{split}
				\aShift_t\left(\cdot,\bar{\eta}_t\right)
				&=
				-\frac{ \Kalman_t\left(\cdot,\bar{\eta}_t\right)\left(H_t+\bar{\eta}_t(H_t)\right)}{2}
				+\frac{\left(
					\left(\Kalman_t\left(\cdot,\bar{\eta}_t\right)\cdot\nabla\right)
					R_t \Kalman_t^{\mathrm{T}}\left(\cdot,\bar{\eta}_t\right)
					\right)^{\mathrm{T}}}{2}
				\\&\phantom{=}+
				\frac{\Kalman_t\left(\cdot,\bar{\eta}_t\right)\mathrm{div}\left(\bar{\eta}_t\tilde{C}_t\Gamma^{\mathrm{T}}_t\right)^{\mathrm{T}}}{2~\bar{\eta}_t}
				+\Omega_t^{0},
			\end{split}
	\end{align}}%
	where $\Omega_t^{0}$ is an arbitrary $\bar{\eta}_t$-harmonic field. Thus the full equation governing the evolution of $\bar{X}$ is given by
	\begin{align}\label{consistent mf equation}
		\begin{split}
			\mathrm{d}
			\bar{X}_t
			&=
			B_t(\bar{X}_t)\mathrm{d}t
			+
			C_t(\bar{X}_t)\mathrm{d}\bar{W}_t
			+
			\tilde{C}_t(\bar{X}_t)\mathrm{d}\bar{V}_t
			\\&\phantom{=}+
			K_t(\bar{X}_t,\bar{\eta}_t)
			\left(\mathrm{d}Y_t
			-\frac{ H_t\left(\bar{X}_t\right)+\int H_t(x)~\bar{\eta}_t(x)~\mathrm{d}x }{2}
			\mathrm{d}t
			\right)
			\\&\phantom{=}+
			\frac{\left(
				\left(K_t\left(\bar{X}_t,\bar{\eta}_t\right)\cdot\nabla\right)
				R_t~
				K_t\left(\bar{X}_t,\bar{\eta}_t\right)^{\mathrm{T}}
				\right)^{\mathrm{T}}
			}{2}\mathrm{d}t
			\\&\phantom{=}+
			\frac{\Kalman_t\left(\bar{X}_t,\bar{\eta}_t\right)\mathrm{div}\left(\bar{\eta}_t\tilde{C}_t\Gamma^{\mathrm{T}}_t\right)^{\mathrm{T}}
				\left(\bar{X}_t\right)
			}{2~\bar{\eta}_t\left(\bar{X}_t\right)}\mathrm{d}t
			+
			\Omega_t^{0}(\bar{X}_t)\mathrm{d}t.
		\end{split}
	\end{align}

	\begin{Rmk}
		We note that
		\begin{align*}
			\frac{\mathrm{div}\left(\bar{\eta}_t\tilde{C}_t\Gamma^{\mathrm{T}}_t\right)}{\bar{\eta}_t}
			&=
			\left(\nabla\log\bar{\eta}_t\right)^{\mathrm{T}}\tilde{C}_t\Gamma^{\mathrm{T}}_t+\mathrm{div}\left(\tilde{C}_t\Gamma^{\mathrm{T}}_t\right)
		\end{align*}
		and thus the correction in \eqref{consistent mf equation} for the correlated observation noise can be rewritten using the gradient of the logarithmic density. 
		This also shows that, just as stated in \cite{NueReiRoz}, if $\tilde{C}_t(x)=\tilde{C}_t$ for all $t\geq 0$, the correlated observation case can be interpreted in terms of the uncorrelated case with a modified observation map $\tilde{H}_t:=H_t-\Gamma_t\tilde{C}^{\mathrm{T}}_t\nabla\log\eta_t$.
		
	\end{Rmk}

	\begin{Rmk}
		Note that in the case of one-dimensional observations $d_y=1$, the term $R_t$ is scalar and $\Kalman_t\left(\cdot,\bar{\eta}_t\right)$ is a $\R^{d_x}$-valued function. One then sees immediately that
		\begin{align*}
			\left(
			\left(\Kalman_t\left(\cdot,\bar{\eta}_t\right)\cdot\nabla\right)
			R_t \Kalman_t^{\mathrm{T}}\left(\cdot,\bar{\eta}_t\right)
			\right)^{\mathrm{T}}
			=
			R_t 
			\left(\Kalman_t\left(\cdot,\bar{\eta}_t\right)\cdot\nabla\right)
			\Kalman_t\left(\cdot,\bar{\eta}_t\right).
		\end{align*}
		
		This is the convective change of the field $\Kalman$ under a flow with velocity $\Kalman$.\newline
		
	\end{Rmk}
	
	\begin{Rmk}
		{We note that for the McKean--Vlasov process $\bar{X}$ in \eqref{consistent mf equation} (and more general any diffusion process driven by $Y$) the tower property of the conditional expectation gives
			\begin{align*}
				\bar{\eta}_t(x):=\mathbb{P}\left(~\bar{X}_t\in\mathrm{d}x~|~Y_{0:T}~\right)
				=
				\mathbb{P}\left(~\bar{X}_t\in\mathrm{d}x~|~Y_{0:t}~\right)~\text{for all}~t\geq 0.
			\end{align*}
			Thus one does not need to fix a timeframe $[0,T]$ a priori and can actually compute $\bar{X}$ online.}
	\end{Rmk}

	\section{The mean field EnKBF in the linear Gaussian setting}\label{EnKBF section}
	
	\begin{Not}
		For any matrix $A$ its symmetric part is denoted by $\mathrm{Sym}\left(A\right):=\frac{A+A^{\mathrm{T}}}{2}$.
	\end{Not}
	
	When $B_t$ and $H_t$ are linear and $C_t,\tilde{C}_t$ are constant matrices for all $t\geq 0$, it is well known that if $\eta_0$ is Gaussian, then the posterior $\eta_t$ will also be Gaussian for all times $t\geq 0$. We denote its (conditional) mean by $\bar{m}_t:=\mathbb{E}_Y\left[\bar{X}_t\right]\in\R^{d_x}$ and its (conditional) covariance by $\bar{P}_t:=\mathbb{Cov}_Y\left[\bar{X}_t,\bar{X}_t\right]\in\R^{d_x\times d_x}$, i.e. $\bar{\eta}_t=\eta_t=\mathcal{N}\left(\bar{m}_t,\bar{P}_t\right)$.\newline
	
	Note that in this case one can easily verify that for any matrix $A\in\R^{d_x\times k}$ with $k\in\N$ arbitrary, it holds that
	\begin{align}\label{divergence of Gaussian}
		\mathrm{div}\left(A\bar{\eta}_t\right)=-\left(x-\bar{m}_t\right)^{\mathrm{T}}\bar{P}_t^{-1}~A \bar{\eta}_t.
	\end{align}
	
	This motivates the Ansatz $\Kalman^0_t=\bar{P}_t H^{\mathrm{T}}_t R^{-1}_{t}$, which in turn results in the Kalman gain
	\begin{align*}
		\Kalman_t=\left(\bar{P}_t H^{\mathrm{T}}_t+\tilde{C}_t\Gamma^{\mathrm{T}}_t\right)R^{-1}_{t}.
	\end{align*}
	
	Since $\Kalman_t$ only depends on the distribution $\bar{\eta}_t$, but not on the state variable $\bar{X}_t$, it follows that
	\begin{align*}
		\left(
		\left(\Kalman_t\left(\cdot,\bar{\eta}_t\right)\cdot\nabla\right)
		R_t \Kalman_t^{\mathrm{T}}\left(\cdot,\bar{\eta}_t\right)
		\right)^{\mathrm{T}}=0
	\end{align*}
	and therefore we obtain by using \eqref{divergence of Gaussian} and by setting $\Omega^0=0$
	\begin{align*}
		\aShift_t\left(x,\bar{\eta}_t\right)
		=
		-\left(\bar{P}_t H^{\mathrm{T}}_t+\tilde{C}_t\Gamma^{\mathrm{T}}_t\right)
		R^{-1}_{t}
		\left(\frac{H_t\left(x+\bar{m}_t\right)}{2}
		+
		\Gamma_t
		\tilde{C}^{\mathrm{T}}_t\bar{P}^{-1}_t\frac{x-\bar{m}_t}{2}
		\right).
	\end{align*}

	{Thus in the linear Gaussian case the following corollary holds.
		
		\begin{Corollary}\label{consistency of EnKBF}
			Assume that $X_0$ is Gaussian, and that for all times $t\geq 0$ both $B_t, H_t$ are linear and both $C_t, \tilde{C}_t$ are constant. Assuming that $\bar{X}$ is the unique solution of the McKean--Vlasov equation
			\begin{align}\label{linear mf ensemble kalman filter}
				\begin{split}
					\mathrm{d}\bar{X}_t
					&=B_t\bar{X}_t\mathrm{d}t
					+C_t\mathrm{d}\bar{W}_t
					+\tilde{C}_t\mathrm{d}\bar{V}_t
					\\&\phantom{=}+
					\left(\bar{P}_t H^{\mathrm{T}}_t+\tilde{C}_t\Gamma^{\mathrm{T}}_t\right)R^{-1}_{t}
					\left(\mathrm{d}Y_t-\frac{H_t\left(\bar{X}_t+\bar{m}_t\right)}{2}\mathrm{d}t\right)
					\\&\phantom{=}-
					\left(\bar{P}_t H^{\mathrm{T}}_t+\tilde{C}_t\Gamma^{\mathrm{T}}_t\right)R^{-1}_{t}
					\Gamma_t\tilde{C}^{\mathrm{T}}_t\bar{P}^{-1}_t\frac{\bar{X}_t-\bar{m}_t}{2}\mathrm{d}t,
				\end{split}
			\end{align}
			where again $\bar{P}_t=\mathbb{Cov}_Y\left[\bar{X}_t\right]$ and $\bar{m}_t=\mathbb{E}_Y\left[\bar{X}_t\right]$,
			satisfying the initial condition $\mathrm{Law}\left(\bar{X}_0\right)=\eta_0=\mathrm{Law}\left(X_0\right)$.\\
			Then it holds that $\mathrm{Law}\left(\bar{X}_t\right)=\eta_t$ for all times $t\geq 0$. Thus $\bar{X}$ is a consistent mean field representation of the posterior $\eta$.
			\begin{flushright}
				$\square$		
			\end{flushright}
		\end{Corollary}
		
		Note that in Corollary \ref{linear mf ensemble kalman filter} we have assumed that equation \eqref{linear mf ensemble kalman filter} is well posed. Indeed since \eqref{linear mf ensemble kalman filter} is a singular McKean--Vlasov equation, that even in the uncorrelated case only satisfies local Lipschitz conditions, this equation falls outside of the standard theory for the well posedness of McKean--Vlasov equations (for example found in \cite{CarDel}). To prove existence and uniqueness of solutions, one can use the relation of \eqref{linear mf ensemble kalman filter} to the famous Kalman--Bucy filter, which is another representation of the posterior in the linear Gaussian setting.\newpage
		
		It is straightforward to derive that the (conditional) mean $\bar{m}$ and the (conditional) covariance $\bar{P}$ of the process $\bar{X}$ in Corollary \ref{consistency of EnKBF} satisfy 
		\begin{subequations}\label{Kalman bucy equations}
			\begin{equation}\label{Kalman Bucy mean}
				\mathrm{d}\bar{m}_t 
				=
				B_t\bar{m}_t\mathrm{d}t+
				\left(\bar{P}_t H^\mathrm{T}_t+\tilde{C}_t\Gamma^{\mathrm{T}}_t\right)R^{-1}_{t}
				\left(\mathrm{d}Y_t-H_t \bar{m}_t\mathrm{d}t\right)
			\end{equation}
			\begin{equation}\label{Kalman Bucy covariance}
				\begin{split}
					\frac{\mathrm{d}\bar{P}_t}{\mathrm{d}t} 
					&=
					B_t\bar{P}_t+\bar{P}_t B^{\mathrm{T}}_t
					+
					C_t C^{\mathrm{T}}_t
					+
					\tilde{C}_t \tilde{C}^{\mathrm{T}}_t
					-
					\left(\bar{P}_t H^\mathrm{T}_t + \tilde{C}_t\Gamma^{\mathrm{T}}_t\right) R^{-1}_t
					\left(H_t \bar{P}_t+\Gamma_t\tilde{C}^{\mathrm{T}}_t\right).
				\end{split}
			\end{equation}
		\end{subequations}
		These are the famous Kalman--Bucy equations \cite[chapter 7, page 228]{Jazwinski} for the correlated noise setting. 
		This again shows the consistency of $\bar{X}$ in the linear Gaussian case. Note that \eqref{Kalman Bucy covariance} is completely decoupled from the observation process and is thus a deterministic quantity.\newline
		
		Of course one can consider system \eqref{Kalman bucy equations} detached from \eqref{linear mf ensemble kalman filter} and easily show the global existence and uniqueness of solutions, which can then be used to prove global existence and uniqueness of solutions to \eqref{linear mf ensemble kalman filter} via a fixed point argument. For the uncorrelated case a well posedness result can also be found in \cite[Remark 2.1]{CriDMJasRuz}, using an argument that was first derived in \cite[Lemma 5.2]{DelMoralTugaut} and can easily be applied to the correlated framework as well. It  proves the following Lemma. For the sake of completeness and the convenience of the reader we also state the proof, arguing just as in \cite{DelMoralTugaut}.\newline
		
		\begin{Lemma}\label{existence prove for linear mf ensemble Kalman}
			Assume that $\bar{P}_0:=\mathbb{Cov}\left[X_0\right]$ is invertible. Then there exists a unique solution $\bar{X}$ to \eqref{linear mf ensemble kalman filter} satisfying the initial condition $\bar{X}_0=X_0$.
		\end{Lemma}
		\begin{proof}
			We make a fixed point argument. Let $(P_t)_{t\geq 0}$ be the solution to \eqref{Kalman Bucy covariance} with initial condition $P_0$. Since the soltions to the matrix Riccati equation \eqref{Kalman Bucy covariance} stay positive definite, $P_t$ is invertible for every time $t\geq 0$.  Now we define the process $\bar{X}^P$ to be the solution to the linear McKean--Vlasov equation
			\begin{align}\label{fixed point equation linear EnKBF}
				\begin{split}
					\mathrm{d}\bar{X}^P_t
					&=B_t~\bar{X}^P_t\mathrm{d}t
					+C_t\mathrm{d}\bar{W}_t
					+\tilde{C}_t\mathrm{d}\bar{V}_t
					\\&\phantom{=}+
					\left(P_t H^{\mathrm{T}}_t+\tilde{C}_t\Gamma^{\mathrm{T}}_t\right)R^{-1}_{t}
					\left(\mathrm{d}Y_t-\frac{H_t\left(\bar{X}^P_t+\mathbb{E}_Y\left[\bar{X}^P_t\right]\right)}{2}\mathrm{d}t\right)
					\\&\phantom{=}-
					\left(P_t H^{\mathrm{T}}_t+\tilde{C}_t\Gamma^{\mathrm{T}}_t\right)
					R^{-1}_{t} \Gamma_t
					\tilde{C}^{\mathrm{T}}_t~P^{-1}_t\frac{\bar{X}^P_t-\mathbb{E}_Y\left[\bar{X}^P_t\right]}{2}\mathrm{d}t
				\end{split}
			\end{align}
			with initial condition $\bar{X}^P_0=X_0$.\newline
			
			We set $\bar{P}_t:=\mathbb{Cov}_Y\left[\bar{X}^P_t,\bar{X}^P_t\right]$. It is easy to see that $\bar{P}$ satsfies the linear equation
			\begin{align*}
				\frac{\mathrm{d}\bar{P}_t}{\mathrm{d}t} 
				&=
				B_t\bar{P}_t+\bar{P}_t B^{\mathrm{T}}_t
				+
				C_t C^{\mathrm{T}}_t
				+
				\tilde{C}_t \tilde{C}^{\mathrm{T}}_t
				-2 
				\mathrm{Sym}\left[
				\left(P_t H^\mathrm{T}_t + \tilde{C}_t\Gamma^{\mathrm{T}}_t\right) R^{-1}_t
				\left(H_t +\Gamma_t\tilde{C}^{\mathrm{T}}_t P^{-1}_t\right)\bar{P}_t\right].
			\end{align*}
			But this equation is also satisfied by $P$, the solution to the Riccati equation \eqref{Kalman Bucy covariance}. Thus by uniqueness of linear equations we derive $P=\bar{P}$ and therefore $\bar{X}^P$ is a solution to \eqref{linear mf ensemble kalman filter}. Given two solutions $\bar{X}^i,~i=1,2$ of \eqref{linear mf ensemble kalman filter}, their (conditional) covariance matrices $\mathbb{Cov}_Y\left[\bar{X}^{i}_t,\bar{X}^{i}_t\right],~i=1,2$ satisfy the Riccati equation \eqref{Kalman Bucy covariance}. Thus by the uniqueness of Riccati equations $\mathbb{Cov}_Y\left[\bar{X}^1_t,\bar{X}^1_t\right]=\mathbb{Cov}_Y\left[\bar{X}^{2}_t,\bar{X}^{2}_t\right]$. Due to the uniqueness of the linear equation \eqref{fixed point equation linear EnKBF}, we can thus conclude $\bar{X}^1=\bar{X}^2$.	
		\end{proof}
		
	}
	
	Unlike for the consistent filter in the general nonlinear setting \eqref{consistent mf equation}, it is clear how \eqref{linear mf ensemble kalman filter} can be approximated using an interacting particle system. For $M\in\N$ let $W^i,~i=1,\cdots,M$ and $V^i,~i=1,\cdots,M$ be independent copies of $W$ and $V$. Let $X^i,~i=1,\cdots,M$ be the solution of
	\begin{align}\label{linear Ensemble Kalman}
		\begin{split}
			\mathrm{d}X^i_t
			&=
			B_t X^i_t\mathrm{d}t+
			C_t\mathrm{d}W^i_t+
			\tilde{C}_t\mathrm{d}V^i_t
			\\&\phantom{=}+
			\left(P^M_t H^\mathrm{T}_t+\tilde{C}_t\Gamma^{\mathrm{T}}_t\right)
			R^{-1}_t
			\left(\mathrm{d}Y_t-\frac{H_t\left(X^i_t+x^M_t\right)}{2}\mathrm{d}t\right)
			\\&\phantom{=}-
			\left(P^M_t H^\mathrm{T}_t+\tilde{C}_t\Gamma^{\mathrm{T}}_t\right)
			R^{-1}_t\Gamma_t
			\tilde{C}^{\mathrm{T}}_t\left(P^M_t\right)^{+}
			\frac{X^i_t-x^M_t}{2}\mathrm{d}t,
		\end{split}
	\end{align}
	where
	\begin{align*}
		x^M_t:=\frac{1}{M}\sum_{i=1}^{M} X^i_t
		~\text{and}~
		P^M_t:=\frac{1}{M-1}\sum_{i=1}^{M}\left(X^i_t-x^M_t\right)\left(X^i_t-x^M_t\right)^\mathrm{T}
	\end{align*}
	denote the ensemble average and the ensemble covariance matrix.  $\left(P^M_t\right)^{+}$ denotes the Moore--Penrose pseudoinverse of $P^M_t$.The pseudoinverse has to be used as for small ensemble sizes $M\leq d_x$, which are commonly used in practice, the ensemble covariance $P^M_t$ can not be invertible.\newline

	\section{The mean field EnKBF for nonlinear signals}\label{EnKBF for nonlinear signals - section}
	For nonlinear signals the posterior can not be expected to be Gaussian and thus the gain term will in general not allow for a explicit description, that admits a simple statistical approximation by interacting particles. Thus deriving suitable numerical approximations to \eqref{consistent mf equation} is challenging. Even though some progress on this subject has been made in recent years (see for example \cite{BerntropGrover} for a Galerkin based approach, \cite{TaghvaeiMehtaMeyn} for a diffusion maps approach and \cite{OlmezTaghvaeiMehta} for a Neural Networks based approach), in practice the EnKBF \eqref{linear Ensemble Kalman}, and its discrete time analogue, is still widely used in the uncorrelated case, even for nonlinear signals, where it will not approximate the optimal filter. A common justification for this is that if the posterior $\eta_t$, and by that also the density of the optimal filter $\bar{\eta}$ are close to a Gaussian, then the mean field EnKBF
	
	\begin{align}\label{nonlinear mf ensemble kalman filter - true inverse}
		\begin{split}
			\mathrm{d}\bar{X}_t
			&=B_t\left(\bar{X}_t\right)\mathrm{d}t
			+C_t\left(\bar{X}_t\right)\mathrm{d}\bar{W}_t
			+\tilde{C}_t\mathrm{d}\bar{V}_t
			\\&\phantom{=}+
			\left(\bar{P}_t H^{\mathrm{T}}_t+\tilde{C}_t\Gamma^{\mathrm{T}}_t\right)R^{-1}_{t}
			\left(\mathrm{d}Y_t-\frac{H_t\left(\bar{X}_t+\bar{m}_t\right)}{2}\mathrm{d}t\right)
			\\&\phantom{=}-
			\left(\bar{P}_t H^{\mathrm{T}}+\tilde{C}_t\Gamma^{\mathrm{T}}_t\right)R^{-1}_{t}
			\Gamma_t\tilde{C}^{\mathrm{T}}_t\bar{P}^{-1}_t\frac{\bar{X}_t-\bar{m}_t}{2}\mathrm{d}t,
		\end{split}
	\end{align}
	where again $\bar{m}_t:=\mathbb{E}_Y\left[\bar{X}_t\right]$ and $\bar{P}_t:=\mathbb{Cov}_Y\left[\bar{X}_t,\bar{X}_t\right]$, will be a good approximation of the optimal filter \eqref{consistent mf equation}.\newline
	
	Indeed if we set $\mathcal{G}\bar{\eta}_t$ to be the normal distribution with mean $\bar{m}_t$ and covariance $\bar{P}_t$, then one obtains \eqref{nonlinear mf ensemble kalman filter - true inverse} by replacing $\bar{\eta}_t$ by $\mathcal{G}\bar{\eta}_t$ in the evolution equation of the optimal filter \eqref{consistent mf equation}. As noted in \cite[Section 2.4.1]{ReichStuart} $\mathcal{G}\bar{\eta}_t$ is indeed the best approximation of $\bar{\eta}$ by a Gaussian with respect to the Kullback--Leibler divergence, i.e.
	\begin{align}\label{Gaussian approximation}
		\mathcal{G}\bar{\eta}_t
		=
		\mathrm{argmin}
		\left\{~
		\int_{\R^{d_x}}~\log\frac{\eta_t}{\mu}~\mathrm{d}\eta_t
		~:~\mu~\text{is Gaussian density on}~\R^{d_x}~
		\right\}.
	\end{align}
	
	Determining conditions when this Gaussian approximation is justified, i.e. when the difference between the optimal filter \eqref{consistent mf equation} and its Gaussian approximation \eqref{nonlinear mf ensemble kalman filter - true inverse} is small, is hard. Even just determining suitable upper bounds of the error based on the Gaussian approximation \eqref{Gaussian approximation} seems difficult as it would involve the Kullback--Leibler.\newline
	
	Therefore we want to discuss another connection between the EnKBF \eqref{nonlinear mf ensemble kalman filter - true inverse} and the consistent filter \eqref{consistent mf equation} that in literature is often referred to as the constant gain approximation \cite{TaghvaeiDeWiljesEtAl}.\newline
	
	\subsection{Constant gain approximation of the optimal filter}\label{subsection constant gain approximation}

	By using partial integration and \eqref{K0} one derives immediately the identity
	\begin{align*}
		\mathbb{E}_{Y}\left[\Kalman^0_t\left(\bar{X}_t,\bar{\eta}_t\right)\right]
		&=
		\int_{\R^{d_x}}\Kalman^0_t\left(x,\bar{\eta}_t\right)~\bar{\eta}_t(x)~\mathrm{d}x
		=
		\int_{\R^{d_x}}\left(x-\bar{m}_t\right)\left(-\mathrm{div}\left(\bar{\eta}_t\Kalman^0_t\left(\cdot,\bar{\eta}_t\right)\right)\right)~\mathrm{d}x
		\\&=
		\int_{\R^{d_x}}\left(x-\bar{m}_t\right)\left(H_t x-H_t \bar{m}_t\right)^{\mathrm{T}}R^{-1}_t\bar{\eta}_t(x)~\mathrm{d}x
		=
		\bar{P}_t H^{\mathrm{T}}_t R^{-1}_t.
	\end{align*}
	
	Thus the the average of $\Kalman$ is given by
	\begin{align}\label{average kalman gain}
		\mathbb{E}_{Y}\left[\Kalman_t\left(\bar{X}_t,\bar{\eta}_t\right)\right]
		=\left(\bar{P}_t H^{\mathrm{T}}_t+\tilde{C}\Gamma^{\mathrm{T}}_t\right) R^{-1}_t.
	\end{align}
	
	\begin{Rmk}
		For nonlinear observations $H_t$ and state dependent $\tilde{C}_t$ the constant gain relation \eqref{average kalman gain} also holds and becomes
		\begin{align*}
			\mathbb{E}_{Y}\left[\Kalman_t\left(\bar{X}_t,\bar{\eta}_t\right)\right]
			=\left(\mathbb{Cov}_{Y}\left[\bar{X}_t,H_t\left(\bar{X}_t\right)\right]+\mathbb{E}_{Y}\left[\tilde{C}_t\left(\bar{X}_t\right)\right]\Gamma^{\mathrm{T}}_t\right) R^{-1}_t.
		\end{align*}
	\end{Rmk}

	While \eqref{average kalman gain} explains the constant gain term in the EnKBF \eqref{nonlinear mf ensemble kalman filter - true inverse}, it does not motivate the approximation of $\nabla\log\bar{\eta}_t(\bar{X}_t)$ by $\bar{P}_t^{-1}\left(\bar{X}_t-\bar{m}_t\right)$. To this end we note that for any affine function $p(x):=\alpha x +\beta$ with $\alpha\in\R^{d_x\times d_x}$ and $\beta\in\R^{d_x}$ one derives by partial integration that
	\begin{align*}
		&\int_{\R^{d_x}} \left(\nabla\log\bar{\eta}_t(x)\cdot p(x)\right)~\bar{\eta}_t(x)~\mathrm{d}x
		=
		\int_{\R^{d_x}} \left(\nabla\bar{\eta}_t(x)\cdot p(x)\right)~\mathrm{d}x
		=
		-\int_{\R^{d_x}} \mathrm{div}(p)~\bar{\eta}_t(x)\mathrm{d}x
		\\&=-\mathrm{tr}\alpha=-\mathrm{tr}\alpha^{\mathrm{T}}
		=-\mathrm{tr}\left[\bar{P}^{-1}_t\bar{P}_t\alpha^{\mathrm{T}}\right]
		=-\mathbb{E}_Y\left[\mathrm{tr}\left[\bar{P}^{-1}_t\left(\bar{X}_t-\bar{m}_t\right)
		\left(\bar{X}_t-\bar{m}_t\right)^{\mathrm{T}}\alpha^{\mathrm{T}}\right]\right]
		\\&=
		-\mathbb{E}_Y\left[\bar{P}^{-1}_t\left(\bar{X}_t-\bar{m}_t\right)
		\cdot\left(\alpha\bar{X}_t-\alpha\bar{m}_t\right)\right]
		=
		\int_{\R^{d_x}} \left(-\bar{P}^{-1}_t(x-\bar{m}_t)\right)\cdot p(x)~\bar{\eta}_t(x)~\mathrm{d}x.
	\end{align*}
	
	Therefore $-\bar{P}^{-1}_t(x-\bar{m}_t)$ is the $\bar{\eta}_t$-orthogonal projection of $\nabla\log\bar{\eta}_t$ onto the space of affine functions.\newline
	
	To summarize, if we denote the $\bar{\eta}_t$-orthogonal projection onto the space of (vector-valued) polynomials of order $k$ by $\pi^{k}[\bar{\eta}_t]$, then
	\begin{align}\label{projection relation}
		\begin{split}
			\pi^{0}[\bar{\eta}_t]\Kalman_t\left(\cdot,\bar{\eta}_t\right)
			&=
			\left(\bar{P}_t H^{\mathrm{T}}_t+\tilde{C}_t\Gamma^{\mathrm{T}}_t\right) R^{-1}_t
			\\
			\pi^{1}[\bar{\eta}_t]\nabla\log\bar{\eta}_t
			&=
			-\bar{P}^{-1}_t(\cdot-\bar{m}_t).
		\end{split}
	\end{align}
	
	Thus, the mean field EnKBF can be seen as an approximation of the optimal filter by projecting its coefficients onto polynomials of certain order. These projections of course cause an inconsistency in the nonlinear setting, however they still result in a similar evolution of the moments. While the projection $\pi^{0}$ will cause a similar evolution of the mean, projections onto the affine functions $\pi^{1}$ will even result in a similar evolution of the covariance. Note that due to the nonlinearity of the signal, this of course does not mean that the mean of the mean field EnKBF is the equal to the one of the optimal filter. However it does explain why the evolution equation of its mean
	\begin{align*}
		\mathrm{d}\bar{m}_t=
		\mathbb{E}_{Y}\left[B_t\left(\bar{X}_t\right)\right]\mathrm{d}t+
		\left(\bar{P}_t H^\mathrm{T}_t+\tilde{C}_t\Gamma^{\mathrm{T}}_t \right) R^{-1}_t
		\left(\mathrm{d}Y_t-H_t \bar{m}_t\mathrm{d}t\right),
	\end{align*}
	is of the same form as the one for the mean of the consistent filter.\newline
	
	\begin{Rmk}
		One of course could use relation \eqref{projection relation} to improve the approximation of the optimal filter by deriving higher order polynomial projections. This question is out of the scope of this paper. In the case of the gain term $\Kalman$ this in particular seems be a non trivial task as it seems to require to look at what happens in \eqref{consistency K} if $\nabla\phi$ is replaced by non-gradient type vector fields. 
	\end{Rmk}
	
	Viewing the mean field EnKBF \eqref{nonlinear mf ensemble kalman filter - true inverse} as a constant gain approximation to the optimal filter \eqref{consistent mf equation} makes it attractive  to quantify its inconsistency by an appropriate coupling of the two differential equation \eqref{nonlinear mf ensemble kalman filter - true inverse} and \eqref{consistent mf equation}. This task is not within the scope of this paper and so we will not investigate this further. However we want to briefly note here that a another more classical way of quantifying the inconsistency is to compare the associated Kolmogorov equation of \eqref{nonlinear mf ensemble kalman filter - true inverse}, which we derive in the next section, with the KSE.
	
	\subsection{Associated stochastic partial differential equation}
	
	Let $\phi$ be an arbitrary smooth test function. Using Itô's rule we derive for the law $\bar{\eta}_t$ of $\bar{X}_t$ defined by equation \eqref{nonlinear mf ensemble kalman filter - true inverse}
	\begin{align}\label{inconsistent PDE}
		\begin{split}
			\mathrm{d}\bar{\eta}_t(\phi)
			&=\bar{\eta}_t(\generator_t\phi)\mathrm{d}t+
			\bar{\eta}_t\left(\nabla\phi\cdot\left(\bar{P}_t H^{\mathrm{T}}_t+\tilde{C}_t\Gamma^{\mathrm{T}}_t\right)\right) R^{-1}_t\mathrm{d}Y_t
			\\&\phantom{=}-
			\bar{\eta}_t\left(\nabla\phi\cdot \left(\bar{P}_t H^{\mathrm{T}}_t+\tilde{C}_t\Gamma^{\mathrm{T}}_t\right) R^{-1}_t \frac{H(\cdot+\bar{m}_t)}{2}\right)\mathrm{d}t
			\\&\phantom{=}-
			\bar{\eta}_t\left(\nabla\phi\cdot \left(\bar{P}_t H^{\mathrm{T}}_t+\tilde{C}_t\Gamma^{\mathrm{T}}_t\right) R^{-1}_t \Gamma_t \tilde{C}^{\mathrm{T}}_t \bar{P}^{-1}_t \frac{\cdot-\bar{m}_t}{2}\right)\mathrm{d}t
			\\&\phantom{=}+
			\frac{1}{2}\bar{\eta}_t\left(\mathrm{tr}\left[\phi'' \left(\bar{P}_t H^{\mathrm{T}}_t+\tilde{C}_t\Gamma^{\mathrm{T}}_t\right) R^{-1}_t \left(H_t\bar{P}_t +\Gamma_t\tilde{C}^{\mathrm{T}}_t\right)\right]\right)
			\mathrm{d}t,
		\end{split}
	\end{align}
	which in general will not coincide with the KSE. In the Gaussian case however we already know that this indeed coincides with the KSE. This can also be shown by using integration by parts on the last four terms in the equation above, as for Gaussian distributions this allows to translate the first and second order differential terms into terms of order zero.\newline

	Before we investigate the well posedness of \eqref{nonlinear mf ensemble kalman filter - true inverse}, we derive some a priori bounds of the (conditional) covariance matrix $\bar{P}$ that will be key for the analysis.\newline

	\subsection{Covariance bounds for the mean field EnKBF}\label{covariance bounds}
	
	To make formulas in this section a bit more concise, we use the following notation.
	\begin{Not}
		For any random vector $Z$ in $\R^{d_x}$  we denote
		\begin{align}\label{definition comCov}
			\begin{split}
				\comCov{B_t}{Z}
				&:=
				\mathbb{E}_{Y}\left[
				\left(B_t(Z)-\mathbb{E}_{Y}[B_t\left(Z\right)]\right)\left(\bar{X}_t-\mathbb{E}_{Y}[Z]\right) ^{\mathrm{T}}\right]
				\\&\phantom{=}+
				\mathbb{E}_{Y}\left[
				\left(\bar{X}_t-\mathbb{E}_{Y}[Z]\right)
				\left(B_t(Z)-\mathbb{E}_{Y}[B_t\left(Z\right)]\right) ^{\mathrm{T}}\right].
			\end{split}
		\end{align}
	\end{Not}
	
	\begin{Not}
		Throughout this paper $\left|\cdot\right|$ shall denote the standard Euclidian norm. If the input is a matrix $A\in\R^{m\times k}$ the result is the Frobeniusnorm $\left| A \right|=\sqrt{\sum_{i=1}^{m}\sum_{j=1}^{k}A_{ij}^2}=\sqrt{\mathrm{tr}\left[AA^{\mathrm{T}}\right]}$. It is easy to see that for any square matrix $A\in\R^{k\times k}$, the trace can be bounded by the Frobeniusnorm in the following way $\mathrm{tr}A\leq \sqrt{k} |A|$. If, on the other hand, $A\in\R^{k\times k}$ is symmetric positive semidefinite, then, by using the singular value decomposition, one can verify easily that $|A|\leq\mathrm{tr}A$.
	\end{Not}

	It is clear that the conditional mean $\bar{m}$ evolves according to
	\begin{align}\label{mean of nonlinear EnKBF}
		\mathrm{d}\bar{m}_t=
		\mathbb{E}_{Y}\left[B_t\left(\bar{X}_t\right)\right]\mathrm{d}t+
		\left(P_t H^\mathrm{T}_t+\tilde{C}_t\Gamma^{\mathrm{T}}_t\right)  R^{-1}_t
		\left(\mathrm{d}Y_t-H_t \bar{m}_t\mathrm{d}t\right).
	\end{align}
	
	Thus Itô's product formula gives us for the centralized process $\bar{X}-\bar{m}$ that
	\begin{align*}
		&\mathrm{d}\left(\bar{X}_t-\bar{m}_t\right)\left(\bar{X}_t-\bar{m}_t\right)^{\mathrm{T}}
		\\&=
		2\mathrm{Sym}\left[
		\left(\bar{X}_t-\bar{m}_t\right) \left(B(\bar{X}_t)-\mathbb{E}_Y\left[B(\bar{X}_t)\right]\right)^{\mathrm{T}}\right]
		\mathrm{d}t
		+
		\left(\left(C(\bar{X}_t)C(\bar{X}_t)^{\mathrm{T}}\right)+\tilde{C}_t\tilde{C}^{\mathrm{T}}\right)\mathrm{d}t
		\\&\phantom{=}
		-
		\mathrm{Sym}
		\left[
		\left(\bar{P}_t H^\mathrm{T}_t+\tilde{C}_t\Gamma^{\mathrm{T}}_t\right)
		R^{-1}_t
		\left( H_t+\Gamma_t\tilde{C}^{\mathrm{T}}_t  \bar{P}^{-1}_t\right)
		\left(\bar{X}_t-\bar{m}_t\right)\left(\bar{X}_t-\bar{m}_t\right)^{\mathrm{T}}
		\right]\mathrm{d}t	
		\\&\phantom{=}+
		2\mathrm{Sym}\left[\left(\bar{X}_t-\bar{m}_t\right)C_t(\bar{X}_t)\mathrm{d}W_t\right]
		+2\mathrm{Sym}\left[\left(\bar{X}_t-\bar{m}_t\right)\tilde{C}_t\mathrm{d}W_t\right]
		.
	\end{align*}
	
	Now by taking the (conditional) expectation $\mathbb{E}_Y$, we see that
	\begin{align}\label{evolution of covariance}
		\begin{split}
			\frac{\mathrm{d}\bar{P}_t}{\mathrm{d}t}
			&=
			\comCov{B_t}{\bar{X}_r}
			+
			\mathbb{E}_{Y}\left[C_t\left(\bar{X}_t\right)C_t\left(\bar{X}_t\right)^{\mathrm{T}}\right]
			+
			\tilde{C}_t \tilde{C}^{\mathrm{T}}_t
			\\&\phantom{=}-
			\left(\bar{P}_t H^\mathrm{T}_t+\tilde{C}_t\Gamma^{\mathrm{T}}_t\right)
			R^{-1}_t
			\left( H_t\bar{P}_t+\Gamma_t\tilde{C}^{\mathrm{T}}_t\right)
		\end{split}
	\end{align}
	and therefore by using that $\mathrm{tr}\left[\left(\bar{P}_t H^\mathrm{T}_t+\tilde{C}_t\Gamma^{\mathrm{T}}_t\right)
	R^{-1}_t
	\left( H_t\bar{P}_t+\Gamma_t\tilde{C}^{\mathrm{T}}_t\right)\right]\geq 0$ one derives
	\begin{align}\label{trace inequality - dynamical}
		\begin{split}
			\frac{\mathrm{d}\mathrm{tr} \bar{P}_t}{\mathrm{d}t}
			&\leq
			2\mathrm{Lip}\left(B\right)
			~\mathrm{tr} \bar{P}_t
			+
			\norm{C_t}{\infty}^2
			+
			\left|\tilde{C}_t\right|^2
			.
		\end{split}
	\end{align}

	Using the deterministic Grönwall Lemma this differential inequality lets us derive	\begin{align}\label{a priori bound covariance}
		\begin{split}
			\sup_{t\leq T} \mathrm{tr} \bar{P}_t
			\leq
			\bar{\Psi}(T)
			&:=
			\exp\left(2~T~\mathrm{Lip}\left(B\right)\right)
			\left(
			\mathrm{tr} \bar{P}_0  
			+	
			\int_{0}^{T} \norm{C_t}{\infty}^2
			+
			\left|\tilde{C}_t\right|^2 ~\mathrm{d}t
			\right).
		\end{split}
	\end{align}
	
	{
		\begin{Rmk}
			We have thus found an upper bound for the spectrum of $\bar{P}$. Note that this is the same as the upper bound for the variance of the signal $X$, i.e. 
			\begin{align*}
				\sup_{t\leq T}\mathrm{tr}\mathbb{Cov}\left[X_t,X_t\right]
				\leq
				\bar{\Psi}(T).
			\end{align*}
			When $\bar{X}$ is a consistent representation of the posterior, i.e. in the linear Gaussian case, this is a stronger version of the inequality implied by the law of total variance. In the inconsistent setting this shows that even though \eqref{nonlinear mf ensemble kalman filter - true inverse} is not the optimal filter, it still satisfies the same conditional variance bound as the true posterior.
		\end{Rmk}
		
		The following positive scalar will play a role in deriving the invertibility of $\bar{P}_t,~t\geq 0$.
		\begin{Def}\label{definition gamma}
			For any $t\in[0,T]$ let
			\begin{align*}
				\gamma_t
				:=
				\lambda_{\min}\left(\tilde{C}_t\left(I-\Gamma^{\mathrm{T}}_t R^{-1}_t \Gamma_t\right)  \tilde{C}^{\mathrm{T}}_t\right)
				+
				\inf_{x}\lambda_{\min}\left(C_t(x) C_t(x)^{\mathrm{T}}\right).
			\end{align*}
			Note that $\lambda_{\min}\left(\tilde{C}_t\left(I-\Gamma^{\mathrm{T}}_t R^{-1}_t \Gamma_t\right)  \tilde{C}^{\mathrm{T}}_t\right)\geq 0$. 
		\end{Def}

		The next Lemma shows that the (conditional) covariance matrix $\bar{P}$ of solutions to \eqref{nonlinear mf ensemble kalman filter - true inverse} stays invertible for all times $t\geq 0$ if $\bar{P}_0$ is already invertible. This hints that the singularity on the right hand side of \eqref{nonlinear mf ensemble kalman filter - true inverse} does not cause a problem, as it will never be seen by the solution $\bar{X}$.\newline}

	\begin{Lemma}\label{regularity of covariance matrix}
		Let $\bar{X}$ be a solution of \eqref{nonlinear mf ensemble kalman filter - true inverse} with conditional covariance matrix $\bar{P}:=\mathbb{Cov}_Y\left[\bar{X},\bar{X}\right]$. Assume that $\bar{P}_0$ is invertible. If the standard Lipschitz conditions in Assumption \ref{standard assumptions} hold and if
		\begin{align}\label{regularity assumption}
			\inf_{t\leq T}\gamma_t
			>0
		\end{align}
		then $\bar{P}_t$ is invertible for all times $t\in[0,T]$. The norm of its inverse can be uniformly bounded by a constant $\underline{\Psi}(T)$, i.e. $\sup_{t\leq T}\left|\left(\bar{P}_t\right)^{-1}\right|\leq\underline{\Psi}(T)$.
	\end{Lemma}
	\begin{proof}
		To not worry too much about the inverse $\bar{P}^{-1}$ in the following computations, we first replace it in equation \eqref{nonlinear mf ensemble kalman filter - true inverse} by the Moore--Penrose pseudoinverse $\bar{P}^+_t$, which is always well defined but may develop singularities. As $\bar{X}$ is assumed to be a solution to \eqref{nonlinear mf ensemble kalman filter - true inverse}, it must also satisfy
		\begin{align}\label{nonlinear mf ensemble kalman filter - pseudo inverse}
			\begin{split}
				\mathrm{d}\bar{X}_t
				&=B_t\left(\bar{X}_t\right)\mathrm{d}t
				+C_t\left(\bar{X}_t\right)\mathrm{d}\bar{W}_t
				+\tilde{C}_t\mathrm{d}\bar{V}_t
				\\&\phantom{=}+
				\left(\bar{P}_t H^{\mathrm{T}}_t+\tilde{C}_t\Gamma^{\mathrm{T}}_t\right)R^{-1}_{t}
				\left(\mathrm{d}Y_t-\frac{H_t\left(\bar{X}_t+\bar{m}_t\right)}{2}\mathrm{d}t\right)
				\\&\phantom{=}-
				\left(\bar{P}_t H^{\mathrm{T}}+\tilde{C}_t\Gamma^{\mathrm{T}}_t\right)R^{-1}_{t}\Gamma_t
				\tilde{C}^{\mathrm{T}}_t\bar{P}^{+}_t\frac{\bar{X}_t-\bar{m}_t}{2}\mathrm{d}t.
			\end{split}
		\end{align}
		
		Since for every $t\geq 0$ the term $\bar{P}_t$ is a symmetric positive semidefinite matrix, there exists an orthogonal matrix $Q_t$ and a diagonal matrix $\Lambda_t:=\mathrm{diag}\left(\lambda^1_t,\cdots,\lambda^{d_x}_t\right)$ such that 
		\begin{align}\label{decomposition P}
			\bar{P}_t=Q^{\mathrm{T}}_t \Lambda_t Q_t~\text{and}~\lambda^i_t\geq 0~\text{for all}~t\geq 0,~i=1,\cdots,d_x.
		\end{align}
		
		It can be shown (see \cite{DieciEirola}) that $Q$ is continuous and $\Lambda_t$ is differentiable and satisfies
		\begin{align}\label{derivative of lambda}
			\frac{\mathrm{d}\Lambda_t}{\mathrm{d}t}
			=
			\mathrm{diag}\left(Q_t\frac{\mathrm{d}\bar{P}_t}{\mathrm{d}t} Q^{\mathrm{T}}_t\right).
		\end{align}
		
		Denote by $e_i$ the $i$-th unit vector in $\R^{d_x}$, then this and \eqref{evolution of covariance} implies
		\begin{align*}
			\frac{\mathrm{d}\lambda^i_t}{\mathrm{d}t}
			&=
			e_i^{\mathrm{T}}
			Q_t\frac{\mathrm{d}\bar{P}_t}{\mathrm{d}t} Q^{\mathrm{T}}_t
			e_i
			\\&=
			e_i^{\mathrm{T}}
			Q_t
			\comCov{B_t}{\bar{X}_t}
			Q^{\mathrm{T}}_t
			e_i
			+
			\mathbb{E}_{Y}
			\left[e_i^{\mathrm{T}}Q_t 	 C_t\left(\bar{X}_t\right)C_t\left(\bar{X}_t\right)^{\mathrm{T}}
			Q^{\mathrm{T}}_te_i
			\right]
			\\&\phantom{=}+
			e_i^{\mathrm{T}}
			Q_t
			\tilde{C}_t \tilde{C}^{\mathrm{T}}_t
			Q^{\mathrm{T}}_t
			e_i
			-
			e_i^{\mathrm{T}}Q_t
			\left(\bar{P}_t H^\mathrm{T}_t+\tilde{C}_t\Gamma^{\mathrm{T}}_t\right)
			R^{-1}_t
			\left( H_t+\Gamma_t\tilde{C}^{\mathrm{T}}_t \bar{P}_t^{+}\right)\bar{P}_t
			Q^{\mathrm{T}}_t
			e_i.
		\end{align*}
		
		We want to estimate $\lambda^i$ from below. To this end we bound all three terms on the right-hand side from below to derive a differential inequality.\newline
		
		First we note that
		\begin{align*}
			\begin{split}
				&e_i^{\mathrm{T}}Q_t \comCov{B_t}{\bar{X}_t} Q^{\mathrm{T}}_te_i
				\\&=
				2~\mathbb{E}_{Y}\left[
				\left(e^{\mathrm{T}}_i Q_t\left(B_t(\bar{X}_t)-B_t(\bar{m}_t)\right)\right)
				\left(e^{\mathrm{T}}_i Q_t\left(\bar{X}_t-\bar{m}_t\right)\right)
				\right]
				\\&\geq
				-2~\sqrt{\mathbb{E}_{Y}\left[
					\left|e^{\mathrm{T}}_iQ^{\mathrm{T}}_t\left(B_t(\bar{X}_t)-B_t(\bar{m}_t)\right) \right|^2
					\right]}
				~\underbrace{\sqrt{\mathbb{E}_{Y}\left[
						\left|e^{\mathrm{T}}_iQ^{\mathrm{T}}_t\left(\bar{X}_t-\bar{m}_t\right)\right|^2
						\right]}}_{=e^{\mathrm{T}}_iQ^{\mathrm{T}}_t\bar{P}_t Q_t e_i=\lambda^{i}_t}
				\\&\geq
				-2\mathrm{Lip}(B)~\sqrt{\mathrm{tr}\bar{P}_t} \sqrt{\lambda^{i}_t}.
			\end{split}
		\end{align*}
		
		By the variational definition of the minimal eigenvalue we also see that
		\begin{align*}
			\begin{split}
				\mathbb{E}_{Y}
				\left[e_i^{\mathrm{T}}Q_t 	 C_t\left(\bar{X}_t\right)C_t\left(\bar{X}_t\right)^{\mathrm{T}}
				Q^{\mathrm{T}}_te_i
				\right]
				&\geq
				\inf_{x}\lambda_{\min}\left(C(x)C(x)^{\mathrm{T}}\right)
				.
			\end{split}
		\end{align*}
		
		Finally it remains to bound
		\begin{align*}
			&e_i^{\mathrm{T}}
			Q_t
			\tilde{C}_t \tilde{C}^{\mathrm{T}}_t
			Q^{\mathrm{T}}_t
			e_i
			-
			e_i^{\mathrm{T}}Q_t
			\left(\bar{P}_t H^\mathrm{T}_t+\tilde{C}_t\Gamma^{\mathrm{T}}_t\right)
			R^{-1}_t
			\left( H_t+\Gamma_t\tilde{C}^{\mathrm{T}}_t \bar{P}_t^{+}\right)\bar{P}_t Q^{\mathrm{T}}_t e_i
			\\&=
			e_i^{\mathrm{T}}Q_t 
			\left(\tilde{C}_t\tilde{C}^{\mathrm{T}}_t-\tilde{C}_t \Gamma^{\mathrm{T}}_t R^{-1}_t \Gamma_t \tilde{C}^{\mathrm{T}}_t \bar{P}_t^{+}\bar{P}_t\right)
			Q^{\mathrm{T}}_t e_i
			-
			e_i^{\mathrm{T}}Q_t\bar{P}_t H^\mathrm{T}_t R^{-1}_t H_t \bar{P}_t Q^{\mathrm{T}}_t e_i
			\\&\phantom{=}-
			e_i^{\mathrm{T}}Q_t
			\left( \tilde{C}_t \Gamma^{\mathrm{T}}_t R^{-1}_t H_t \bar{P}_t 
			+
			\bar{P}_t H^\mathrm{T}_t R^{-1}_t \Gamma_t \tilde{C}^{\mathrm{T}}_t \bar{P}_t^{+}\bar{P}_t
			\right)
			Q^{\mathrm{T}}_t
			e_i
		\end{align*}
		from below.	The decomposition \eqref{decomposition P} of $\bar{P}_t$, lets us also easily derive 
		\begin{align*}
			e_i^{\mathrm{T}}Q_t\bar{P}_t H^\mathrm{T}_t R^{-1}_t H_t \bar{P}_t Q^{\mathrm{T}}_t e_i
			&\leq \lambda_{\max}\left(H^\mathrm{T}_t R^{-1}_t H_t\right) \left(\lambda^{i}_t\right)^2 
			\\
			e_i^{\mathrm{T}}Q_t
			\left( \tilde{C}_t \Gamma^{\mathrm{T}}_t R^{-1}_t H_t \bar{P}_t 
			+
			\bar{P}_t H^\mathrm{T}_t R^{-1}_t \Gamma_t \tilde{C}^{\mathrm{T}}_t \bar{P}_t^{+}\bar{P}_t
			\right)
			Q^{\mathrm{T}}_t
			e_i
			&\leq
			2\left|\tilde{C}_t \Gamma^{\mathrm{T}}_t R^{-1}_t H_t\right| \lambda^i_t
			\\
			e_i^{\mathrm{T}}Q_t 
			\left(\tilde{C}_t\tilde{C}^{\mathrm{T}}_t-\tilde{C}_t \Gamma^{\mathrm{T}}_t R^{-1}_t \Gamma_t \tilde{C}^{\mathrm{T}}_t \bar{P}_t^{+}\bar{P}_t\right)
			Q^{\mathrm{T}}_t e_i
			&\geq
			\lambda_{\min}\left(\tilde{C}_t \left(I-\Gamma^{\mathrm{T}}_t R^{-1}_t \Gamma_t\right) \tilde{C}^{\mathrm{T}}_t\right).
		\end{align*}
		Note that for the last inequality we also used that $\tilde{C}_t \Gamma^{\mathrm{T}}_t R^{-1}_t \Gamma_t \tilde{C}^{\mathrm{T}}_t$  is positive semidefinite. And thus we derive
		\begin{align}\label{lower bound lambda}
			\begin{split}
				\frac{\mathrm{d}\lambda^i_t}{\mathrm{d}t}
				&\geq
				-2\mathrm{Lip}(B)~\sqrt{\mathrm{tr}\bar{P}_t} \sqrt{\lambda^{i}_t}
				\\&\phantom{=}+
				\underbrace{\inf_{x}\lambda_{\min}\left(C_t(x)C_t(x)^{\mathrm{T}}\right)
					+
					\lambda_{\min}\left(\tilde{C}_t \left(I- \Gamma^{\mathrm{T}}_t R^{-1}_t \Gamma_t\right) \tilde{C}^{\mathrm{T}}_t\right)}_{\geq\inf_{s\leq T}\gamma_s}
				\\&\phantom{=}-
				\lambda_{\max}\left(H^\mathrm{T}_t R^{-1}_t H_t\right) \left(\lambda^{i}_t\right)^2 
				-2
				\left|\tilde{C}_t \Gamma^{\mathrm{T}}_t R^{-1}_t H_t\right|
				\lambda^i_t
				.
			\end{split}
		\end{align}
		
		Let $\underline{\lambda}$ denote a solution of
		\begin{align}\label{lower bound lambda 2}
			\begin{split}
				\frac{\mathrm{d}\underline{\lambda}_t}{\mathrm{d}t}
				&=
				-2\mathrm{Lip}(B)~\sqrt{\mathrm{tr}\bar{P}_t} \sqrt{\underline{\lambda}_t}
				+
				\frac{\inf_{s\leq T}\gamma_s}{2}
				\\&\phantom{=}-
				\lambda_{\max}\left(H^\mathrm{T}_t R^{-1}_t H_t\right) \left(\underline{\lambda}_t\right)^2 
				-2
				\left|\tilde{C}_t \Gamma^{\mathrm{T}}_t R^{-1}_t H_t\right|
				\underline{\lambda}_t
				.
			\end{split}
		\end{align}
		with initial condition $\underline{\lambda}_0=\frac{\min_{i}\lambda^i_0}{2}$. The existence of solutions (which might not be unique) is guaranteed by the Peano theorem. 
		Then by \cite[II. Lemma, page 64]{Walter}, $\underline{\lambda}$ is a lower bound for all $\lambda^i$.
		
		By assumption \eqref{regularity assumption}, the right-hand side of \eqref{lower bound lambda 2} is strictly bigger than zero, when $\underline{\lambda}_t=0$. Thus for $t_0>0$ with $\underline{\lambda}_{t_0}=0$, we have $\underline{\lambda}'(t_0)>0$. According to \cite[II. Lemma, page 64]{Walter}, this indeed implies $\underline{\lambda}_t>0$ for all $t\geq 0$, if $\underline{\lambda}_0>0$, which is the case if $\bar{P}_0$ is regular.
		Since $\left|\bar{P}^{-1}_t\right|^2=\sum_{i=1}^{d_x}\left(1/\lambda^i_t\right)^2$ this indeed gives us a bound $\underline{\Psi}(T):=\sqrt{d_x}\frac{1}{\inf_{t\leq T} \underline{\lambda}_t}$ for the inverse of the covariance.
	\end{proof}
	
	\begin{Rmk}
		Note that if $\bar{P}_0$ is not invertible, then one can still show that $\bar{P}_t$ is invertible for all $t>0$. To accomplish this one uses the proof above and argues by contradiction. As \cite[II. Lemma, page 64]{Walter} also implies that, if there existed a $t_0>0$ with $\underline{\lambda}_{t_0}=0$, then there would be a $\bar{t}>0$ such that $\underline{\lambda}_t\leq 0$ for all $t\leq\bar{t}$. We know that $\lambda^i$ is non-negative, thus it would hold that $\underline{\lambda}_t= 0$ for all $t\leq\bar{t}$.\\
		However by \eqref{lower bound lambda 2} and assumption \eqref{regularity assumption} this would also mean that $\frac{\mathrm{d}\underline{\lambda}_t}{\mathrm{d}t}>0$ on $\left(0,\bar{t}\right)$, which is a contradiction to $\underline{\lambda}$ being constant on $\left(0,\bar{t}\right)$. Therefore we indeed can conclude that $\underline{\lambda}_t>0$ and thus $\lambda^i_t>0$  for all $t> 0,~i=1,\cdots,d_x$. Therefore $\bar{P}_t$ is regular.\newline
	\end{Rmk}

	Next we investigate the well posedness of \eqref{nonlinear mf ensemble kalman filter - true inverse}.
	
	\subsection{Well posedness of the mean field Ensemble Kalman--Bucy filter}\label{section proof of thm}

	One difficulty in the analysis of \eqref{nonlinear mf ensemble kalman filter - true inverse} that is apparent at first sight is, that this equation is (potentially) singular, i.e. the inverse of the (conditional) covariance on the right hand side is not well defined for all probability distributions and blows up when $\bar{P}_{t_0}$ becomes non invertible for some point in time $t_0$.\newline
	
	However even in the uncorrelated setting $\tilde{C}=0$ we note that the coefficients in \eqref{nonlinear mf ensemble kalman filter - true inverse} are only locally Lipschitz and thus this equation falls outside the standard theory of McKean--Vlasov equations that is e.g. found in \cite{CarDel}. For ordinary SDEs with local Lipschitz conditions the main tool for proving well posedness is the introduction of stopping times that localize the equation and reduce its analysis to the case of global Lipschitz conditions.	However looking at a stopped stochastic process of course changes its time marginal laws and due to the law dependence of McKean--Vlasov equations such stopping times would thus change the dynamics of the equation. So far there does not seem to exist a localization argument that is suited to prove the well posedness of a general class of McKean--Vlasov equations under local Lipschitz assumptions. Indeed \cite{Scheutzow2} shows that such an argument can not hold in full generality as it constructs locally Lipschitz McKean--Vlasov differential equations (without the influence of noise) where uniqueness does not hold. Nevertheless in recent years some progress has been made regarding the well posedness of McKean--Vlasov equations under local Lipschitz conditions. For example in \cite{HammersleySiskaSzpruch} existence of weak solutions has been shown by the usage of Lyapunov techniques. \cite{Liu} show the existence of solutions for a fairly general class of equations and determine some conditions for their uniqueness, but even in the uncorrelated case these results do not apply to \eqref{nonlinear mf ensemble kalman filter - true inverse},  in parts due to missing growth conditions. \newline
	
	Before we show our main well posedness result in the correlated noise framework, we first look at the uncorrelated case and make some restricting assumptions regarding the observations, as it will motivate the partial stopping argument used later on in the proof for the general case.\newline
	
	\begin{Lemma}\label{Existence and uniqueness for mf - special case}
		Beside the standard Lipschitz conditions (see assumption \ref{standard assumptions}), we assume that $C$ is bounded and that all coefficients are continuously differentiable with respect to $t$. Furthermore we assume that $\tilde{C}=0$ and that there exists some scalar valued function $\alpha$ such that $H^\mathrm{T}_t R^{-1}_t H_t=\alpha_t I$ for all $t\geq 0$, where $I$ again denotes the identity matrix. Then there exists a unique solution $\bar{X}$ of \eqref{nonlinear mf ensemble kalman filter - true inverse} on the time interval $[0,T]$.
	\end{Lemma}

	To prove Lemma \ref{Existence and uniqueness for mf - special case} we will make a fixed point argument just as we did in the linear case. However, since the equation of the (conditional) covariance matrix $\bar{P}$ does not decouple from $\bar{X}$ for nonlinear $B$, we will not be able to simply guess the right fixed point. Instead we have to work with a priori bounds for the first two moments which we derive now, before proving Lemma \ref{Existence and uniqueness for mf - special case}.\newline

	{First we define the set of suitable fixed points.
		
		\begin{Def}\label{definition spsd}
			We denote by $\mathrm{SPSD}^1_Y([0,T])$ the set of all differentiable processes $P$ such that
			\begin{itemize}
				\item $P$ is adapted to the natural filtration of $Y$,
				
				\item  $P_t$ is a symmetric positive semidefinite matrix for every $t\geq 0$,
				
				\item $P_0=\bar{P}_0$.
			\end{itemize} 
		\end{Def}		
		
	} For any $P\in\mathrm{SPSD}^1_Y([0,T])$ we define the stochastic process $\bar{X}^P$ to be the solution of
	\begin{align*}
		\mathrm{d}\bar{X}^P_t
		&=
		B_t(\bar{X}^P_t)\mathrm{d}t+
		C_t(\bar{X}^P_t)\mathrm{d}\bar{W}_t
		+
		P_t H^\mathrm{T}_t
		R^{-1}_t
		\left(\mathrm{d}Y_t-\frac{H_t\left(\bar{X}^P_t+\bar{m}^P_t\right)}{2}\mathrm{d}t
		\right).
	\end{align*}
	
	Hereby we denote $\bar{m}^P_t:=\mathbb{E}_{Y}\left[\bar{X}^P_t\right]$. Since this is a McKean--Vlasov equation satisfying the usual Lipschitz conditions, there indeed exists a unique strong solution \cite{CarDel}.\newline
	
	A solution of \eqref{nonlinear mf ensemble kalman filter - true inverse} is given by a fixed point, i.e. a suitable process $P$ such that
	\begin{align*}
		\mathbb{Cov}_{Y}\left[\bar{X}^P_t\right]=P_t~\text{for all}~t\in[0,T].
	\end{align*}
	We will prove its existence and uniqueness on small time intervals with the Banach fixed point theorem. The extension to arbitrary time intervals then follows from a standard glueing argument. In order to prove that the map 
	\begin{align*}
		P\mapsto\mathbb{Cov}_{Y}\left[\bar{X}^P\right]
	\end{align*}
	is a strict contraction on a suitable subset of $\mathrm{SPSD}^1_Y([0,T])$, a priori estimates on the (conditional) covariance matrix $\mathbb{Cov}_{Y}\left[\bar{X}^P_t\right]$, independent of $P$ and $t$, will be key.\newline
	
	We note that the (conditional) mean and covariance matrix of $\bar{X}^P$ satisfy the following two equations
	\begin{align}\label{mean of XP}
		\mathrm{d}\bar{m}^P_t=
		\mathbb{E}_{Y}\left[B_t\left(\bar{X}^P_t\right)\right]\mathrm{d}t+
		P_t H^\mathrm{T}_t  R^{-1}_t
		\left(\mathrm{d}Y_t-H_t \bar{m}^P_t\mathrm{d}t\right),
	\end{align}
	and
	\begin{align}\label{covariance of XP}
		\begin{split}
			\frac{\mathrm{d}\mathbb{Cov}_{Y}\left[\bar{X}^P_t\right]}{\mathrm{d}t}
			&=
			\comCov{B_t}{\bar{X}^P_t}
			+
			\mathbb{E}_{Y}\left[C_t\left(\bar{X}^P_t\right)C_t\left(\bar{X}^P_t\right)^{\mathrm{T}}\right]
			\\&\phantom{=}-
			\mathrm{Sym}\left(
			P_t H^\mathrm{T}_t
			R^{-1}_t
			H_t\mathbb{Cov}_{Y}\left[\bar{X}^P_t\right]
			\right).
		\end{split}
	\end{align}
	
	\begin{Rmk}\label{robustness as a lyapunov function}
		We note that if $H^\mathrm{T}_t R^{-1}_t H_t$ is not a scalar matrix, that is there exists no scalar-valued function $\alpha$ such that $H^\mathrm{T}_t R^{-1}_t H_t =\alpha_t I$, then we can not guarantee that $P_t H^\mathrm{T}_t R^{-1}_t H_t$ is positive definite, let alone symmetric. Thus we can not simply take the trace in \eqref{covariance of XP}, use the Grönwall lemma and expect to derive bounds independent of $P$. In other words, the covariance matrix is a Lyapunov function of \eqref{nonlinear mf ensemble kalman filter - true inverse}, that, for general coefficients, is not robust to perturbations in the covariance matrix.
	\end{Rmk}

	When the assumptions of Lemma \ref{Existence and uniqueness for mf - special case}  hold, equation \eqref{covariance of XP} reduces to
	\begin{align}\label{covariance of XP simple case}
		\begin{split}
			\frac{\mathrm{d}\mathbb{Cov}_{Y}\left[\bar{X}^P_t\right]}{\mathrm{d}t}
			&=
			\comCov{B_t}{\bar{X}^P_t}+
			\mathbb{E}_{Y}\left[C_t\left(\bar{X}^P_t\right)C_t\left(\bar{X}^P_t\right)^{\mathrm{T}}\right]
			\\&-
			\frac{\alpha_t}{2} 
			\left(P_t\mathbb{Cov}_{Y}\left[\bar{X}^P_t\right]
			+
			\mathbb{Cov}_{Y}\left[\bar{X}^P_t\right] P_t\right).
		\end{split}
	\end{align}
	We note that by the cyclical invariance of the trace we have
	\begin{align*}
		\mathrm{tr}\left[P_t\mathbb{Cov}_{Y}\left[\bar{X}^P_t\right]\right]
		&=
		\mathbb{E}_{Y}\left[
		\mathrm{tr}\left[P_t \left(X^P_t-m^P_t\right)\left(X^P_t-m^P_t\right)^{\mathrm{T}}\right]
		\right]
		\\&=
		\mathbb{E}_{Y}\left[
		\left(X^P_t-m^P_t\right)^{\mathrm{T}} P_t \left(X^P_t-m^P_t\right)
		\right]
		\geq 0,
	\end{align*}
	which gives us  the differential inequality
	\begin{align*}
		\frac{\mathrm{d}\mathrm{tr}\mathbb{Cov}_{Y}\left[\bar{X}^P_t\right]}{\mathrm{d}t}
		\leq
		2~\mathrm{Lip}(B_t)~\mathrm{tr}\mathbb{Cov}_{Y}\left[\bar{X}^P_t\right] +\norm{C_t}{\infty}^2,
	\end{align*}
	and thus we derive by the deterministic Grönwall inequality
	\begin{align}\label{bound for cov trace - local}
		\sup_{t\leq T}~\mathrm{tr}~\mathbb{Cov}_{Y}	\left[X^P_t\right]
		\leq
		\exp\left(2T~\mathrm{Lip}(B)\right) \left(\mathrm{tr} \bar{P}_0+\int_{0}^{T}\norm{C_t}{\infty}^2\mathrm{d}t\right)
		=:\kappa_0(T)
		.
	\end{align}
	
	This means that we can assume that $\sup_{t\leq T} \mathrm{tr}P_t\leq\kappa_0(T)$. Using \eqref{covariance of XP simple case} we thus obtain
	\begin{align}\label{bound derivative of covariance}
		\sup_{t\leq T}\left| \frac{\mathrm{d}\mathrm{tr}\mathbb{Cov}_{Y}\left[X^P_t\right]}{\mathrm{d}t}\right|
		\leq
		\sup_{t\leq T}\left(2\mathrm{Lip}(B)+\alpha_t\right) \kappa_0(T)^2=:\kappa_1(T).
	\end{align}
	
	With the two bounds \eqref{bound for cov trace - local} and \eqref{bound derivative of covariance}, we define the set
	\begin{align}\label{set of potential fixed points}
		\begin{split}
			\mathcal{X}_T:=\left\{~P\in\mathrm{SPSD}^1_Y([0,T])
			~:~
			\norm{P}{\infty}\leq\kappa_0(T),~\norm{\partial_t P}{\infty}\leq\kappa_1(T)
			\right\}.
		\end{split}
	\end{align}
	
	The following Lemma will be useful for deriving a contraction property for the fixed point map $P\mapsto\mathbb{Cov}_Y\left[\bar{X}^P,\bar{X}^P\right]$.
	\begin{Lemma}\label{bound rough integral difference}
		For every fixed $T>0$ there exists a constant $\kappa_{Y}(T)$, such that 
		\begin{align*}
			\sup_{t\leq T}\left|\int_{0}^{t} (P^1_s-P^2_s) H^\mathrm{T}_s R^{-1}_s\mathrm{d}Y_s\right|
			\leq \kappa_{Y}(T) \norm{P^1-P^2}{\Cspace^1}
		\end{align*}
		holds for all $P^1,P^2\in\mathrm{SPSD}^1_Y([0,T])$. Furthermore these constants satisfy $\kappa_{Y}(T)\xrightarrow{T\to 0} 0$.
	\end{Lemma}
	\begin{proof}
		We note that all $P^1,P^2\in\mathrm{SPSD}^1_Y([0,T])$ are processes of bounded variation and thus one gets by integration by parts
		\begin{align*}
			&\int_{0}^{t}  (P^1_s-P^2_s)H^\mathrm{T}_s R^{-1}_s \mathrm{d}Y_s
			=
			\left(P^1_t-P^1_t\right) H^\mathrm{T}_t R^{-1}_t	Y_t
			-
			\int_{0}^{t} \frac{\mathrm{d}~(P^1_s-P^1_2)H^\mathrm{T}_s R^{-1}_s}{\mathrm{d}s} Y_s\mathrm{d}s,
		\end{align*}
		where the boundary term at $t=0$ disappears as $Y_0=0$.\newline
		
		Since $H^\mathrm{T}_sR^{-1}_s$ is a $\Cspace^1\left([0,T];\R^{d_x\times d_y}\right)$ function we immediately derive the desired inequality with
		\begin{align*}
			\kappa_{Y}(T) 
			:=
			\sup_{t\leq T}\left|H^\mathrm{T}_t R^{-1}_t	Y_t\right|
			+
			\int_{0}^{T} \left|H^\mathrm{T}_t R^{-1}_t	Y_t\right| +
			\left|\frac{\mathrm{d}H^\mathrm{T}_t R^{-1}_t}{\mathrm{d}t}	Y_t\right|~\mathrm{d}t.
		\end{align*} 
		
		By the continuity of $Y$ and $Y_0=0$ we also immediately see that $\kappa_{Y}(T)\xrightarrow{T\to 0} 0$.
	\end{proof}
	
	As a corollary one also immediately derives for any $P\in\mathcal{X}_T$ that
	\begin{align}\label{bound rough integral}
		\sup_{t\leq T}\left|\int_{0}^{t} P_s H^\mathrm{T}_s R^{-1}_s\mathrm{d}Y_s\right|
		\leq \kappa_{Y}(T) \sqrt{\kappa_0(T)^2+\kappa_1(T)^2}.
	\end{align}
	
	Using the Lipschitz continuity of $B$ we see that for any $P\in\mathcal{X}_T$
	\begin{align*}
		&\left|m^P_t\right|
		\leq
		\left|m^P_0\right|
		+
		\int_{0}^{t} \left|B_s(0)\right|\mathrm{d}s+
		\int_{0}^{t}
		\mathrm{Lip}(B)
		\mathbb{E}_Y\left[\left|X^P_s-m^P_s\right|\right]\mathrm{d}s
		\\&\phantom{=}+
		\left|\int_{0}^{t}P_s H^\mathrm{T}_s R^{-1}_s\mathrm{d}Y_s\right|
		+
		\int_{0}^{t} \big(\underbrace{\left|P_s\right|}_{\leq\kappa_0(T)} \left|H^\mathrm{T}_s R^{-1}_s H_s\right|+\mathrm{Lip}(B)\big)
		\left|m^P_s\right|\mathrm{d}s.
	\end{align*}
	
	Applying the deterministic Grönwall inequality together with
	\begin{align*}
		\mathbb{E}_{Y}\left[\left|X^P_s-m^P_s\right|\right]
		\leq
		\sqrt{\mathbb{E}_{Y}\left[\left|X^P_s-m^P_s\right|^2\right]}
		=\sqrt{\mathrm{tr}~\mathbb{Cov}_{Y}\left[X^P_s\right]}
		\leq \sqrt{\kappa_0(T)}
	\end{align*}
	and \eqref{bound rough integral}, we obtain for $P\in\mathcal{X}_T$
	\begin{align*}
		\sup_{t\leq T}\left|m^P_t\right|
		&\leq\kappa_{m}(T)
		:=
		\exp\left(
		\int_{0}^{t} \left( \kappa_0(T) \left|H^\mathrm{T}_s R^{-1}_s H_s\right|+\mathrm{Lip}(B)\right)
		\mathrm{d}s\right)
		\\&
		\phantom{=}\times
		\left(
		\left|m^P_0\right|
		+
		\int_{0}^{t}
		\left|B_s(0)\right|+
		\mathrm{Lip}(B)
		\sqrt{\kappa_0(T)}\mathrm{d}s+
		\kappa_{Y}(T) \sqrt{\kappa_0(T)^2+\kappa_1(T)^2}
		\right).
	\end{align*}

	Now we are in the position to prove Lemma \ref{Existence and uniqueness for mf - special case}.
	\begin{proof}[Proof of Lemma \ref{Existence and uniqueness for mf - special case}]
		As already noted we look for a fixed point $\bar{P}$ such that $\bar{P}=\mathbb{Cov}_Y\left[\bar{X}^{\bar{P}},\bar{X}^{\bar{P}}\right]$. Due to the previously derived bounds \eqref{bound for cov trace - local} and \eqref{bound derivative of covariance}, we can restrict our search to the set $\mathcal{X}_T$.	As a closed subset of a Banach space, $\mathcal{X}_T$ is itself a complete metric space.\newline
		
		Thus we can employ the Banach fixed point theorem and it only remains to show the contraction property of the fixed point map $P\mapsto\mathbb{Cov}_Y\left[\bar{X}^P,\bar{X}^P\right]$ on $\mathcal{X}_T$ for sufficiently small $T$. To this end let $P^1, P^2\in\mathcal{X}_T$ be given and denote by $X^1, X^2$ the corresponding solutions and by $m^1, m^2$ their means. Then it is clear that
		\begin{align*}
			X^1_t-X^2_t
			&=
			\int_{0}^{t}\left(B_s(X^1_s)-B_s(X^2_s)\right)\mathrm{d}s
			+
			\int_{0}^{t}\left(P^1_s-P^2_s\right)H^{\mathrm{T}}_s R^{-1}_s\mathrm{d}Y_s
			\\&\phantom{=}-
			\int_{0}^{t}\left(P^1_s-P^2_s\right)H^{\mathrm{T}}_s R^{-1}_s H_s \frac{X^1_s+m^1_s}{2}
			\mathrm{d}s
			\\&\phantom{=}-
			\int_{0}^{t}P^2_s H^{\mathrm{T}}_s R^{-1}_s H_s \frac{X^1_s-X^2_s+m^1_s-m^2_s}{2}
			\mathrm{d}s
			\\&\phantom{=}+
			\int_{0}^{t}\left(C_s(X^1_s)-C_s(X^2_s)\right)\mathrm{d}\bar{W}_s
		\end{align*}

		Using the estimates $\left|\sum_{i=1}^{k} a_i\right|^2 ~\leq~k~  \sum_{i=1}^{k} \left|a_i\right|^2$ and $\left|\int_{0}^{t} f(t)\mathrm{d}t\right|^2\leq t\int_{0}^{t} \left|f(t)\right|^2\mathrm{d}t$, we get 
		\begin{align*}
			\left|X^1_t-X^2_t\right|^2
			&\leq
			5t\int_{0}^{t}\left|B_s(X^1_s)-B_s(X^2_s)\right|^2\mathrm{d}s
			+
			5
			\left|\int_{0}^{t}\left(P^1_s-P^2_s\right)H^{\mathrm{T}}_s R^{-1}_s \mathrm{d}Y_s\right|^2
			\\&\phantom{=}+
			5 t\int_{0}^{t}\left|P^1_s-P^2_s\right|^2\left|H^{\mathrm{T}}_s R^{-1}_s H_s \right|^2\frac{|X^1_s|^2+|m^1_s|^2}{2}
			\mathrm{d}s
			\\&\phantom{=}+
			5t\int_{0}^{t}\left|P^2_s\right|^2 \left|H^{\mathrm{T}}_s R^{-1}_s H_s\right|^2
			\frac{\left|X^1_s-X^2_s\right|^2+\left|m^1_s-m^2_s\right|^2}{2}
			\mathrm{d}s
			\\&\phantom{=}+
			5\left|\int_{0}^{t}\left(C_s(X^1_s)-C_s(X^2_s)\right)\mathrm{d}\bar{W}_s\right|^2
		\end{align*}
		
		Note that
		\begin{align*}
			\left|m^1_s-m^2_s\right|^2 \leq \mathbb{E}_{Y}\left[\left|X^1_s-X^2_s\right|^2\right]
			~\text{and}~
			\mathbb{E}_{Y}\left[|X^1_s|^2\right]+|m^1_s|^2
			\leq 
			\kappa_0(T)+2\kappa_{m}(T)^2.
		\end{align*}
		
		Using the Lipschitz continuity of $B$ and $C$, as well as Itô isometry, we thus derive
		\begin{align*}
			&\mathbb{E}_{Y}\left[\left|X^1_t-X^2_t\right|^2\right]	
			\\&\leq
			5\int_{0}^{t}
			\left(T\mathrm{Lip}(B)^2+T\left|P^2_s\right|^2\left|H^{\mathrm{T}}_s R^{-1}_s H_s\right|^2
			+\mathrm{Lip}(C)^2\right)
			~\mathbb{E}_{Y}\left[\left|X^1_s-X^2_s\right|^2\right]	\mathrm{d}s
			\\&\phantom{=}+
			5
			\left(\kappa_{Y}(T)^2 + 5 T \frac{\kappa_0(T)+2\kappa_{m}(T)^2}{2}
			\int_{0}^{t}\left|H^{\mathrm{T}}_s R^{-1}_s H_s\right|^2\mathrm{d}s\right)
			\norm{P^1-P^2}{\Cspace^1}^2.
		\end{align*}
		
		Thus we derive from the Grönwall lemma
		\begin{align*}
			&\sup_{t\leq T}
			\mathbb{E}_{Y}\left[\left|X^1_t-X^2_t\right|^2\right]
			\\&\leq
			\exp\left(
			5\int_{0}^{T}
			T\mathrm{Lip}(B)^2+T\left|P^2_s\right|^2\left|H^{\mathrm{T}}_s R^{-1}_s H_s\right|^2
			+\mathrm{Lip}(C)^2\mathrm{d}s
			\right)
			\\&\phantom{=}\times
			5
			\left(\kappa_{Y}(T)^2 + 5 T \frac{\kappa_0(T)+2\kappa_{m}(T)^2}{2}
			\int_{0}^{T}\left|H^{\mathrm{T}}_s R^{-1}_s H_s\right|^2\mathrm{d}s\right)
			\norm{P^1-P^2}{\Cspace^1}^2
			\\&:=\kappa_{\mathrm{contr}}(T)\norm{P^1-P^2}{\Cspace^1}^2.
		\end{align*}
		
		Where clearly $\kappa_{\mathrm{contr}}(T)\xrightarrow{T\to 0} 0$. By employing 
		\begin{align*}
			&\left|
			\mathbb{Cov}_{Y}\left[X^1_t\right]
			-
			\mathbb{Cov}_{Y}\left[X^2_t\right]
			\right|
			\leq
			2\sqrt{
				\mathrm{tr}\mathbb{Cov}_{Y}\left[X^1_t\right]
				+
				\mathrm{tr}\mathbb{Cov}_{Y}\left[X^2_t\right]
			}
			\sqrt{\mathbb{E}_{Y}\left[\left|X^1_t-X^2_t\right|^2\right]}
		\end{align*}
		we therefore obtain
		\begin{align*}
			\sup_{t\leq T}\left|
			\mathbb{Cov}_{Y}\left[X^1_t\right]
			-
			\mathbb{Cov}_{Y}\left[X^2_t\right]
			\right|
			\leq
			2\sqrt{2\kappa_0(T)\kappa_{\mathrm{contr}}(T)}\norm{P^1-P^2}{\Cspace^1},
		\end{align*}
		
		Using the differential equation for the evolution of the covariance matrix \eqref{covariance of XP simple case} and the previously derived bounds, one also derives the existence of $q(T)>0$ with $q(T)\xrightarrow{T\to 0} 0$ such that
		\begin{align*}
			\sup_{t\leq T}\left|
			\frac{\mathrm{d}\mathbb{Cov}_{Y}\left[X^1_t\right]}{\mathrm{d}t}
			-
			\frac{\mathrm{d}\mathbb{Cov}_{Y}\left[X^2_t\right]}{\mathrm{d}t}
			\right|
			\leq
			q(T)\norm{P^1-P^2}{\Cspace^1}.
		\end{align*}
		
		Thus we have proven the contraction property for sufficiently small $T>0$ and therefore Lemma \ref{Existence and uniqueness for mf - special case} holds for small time domains. To prove this fact for arbitrary timeframes, one now simply uses a standard glueing argument.\newline
		
	\end{proof}

	Next we generalize Lemma \ref{Existence and uniqueness for mf - special case} to the correlated noise framework without assuming that $H^\mathrm{T}_t R^{-1}_t H_t$ is scalar. To this end we have to introduce a partial stopping argument that makes the trace of the covariance robust as a Lyapunov function (see Remark \ref{robustness as a lyapunov function}), since this was the main argument that allowed us to use the Banach fixed point theorem. We now are in the position to show our most general well posedness result for \eqref{nonlinear mf ensemble kalman filter - true inverse}.

	\begin{Thm}\label{Existence and uniqueness for mf}
		Beside the standard Lipschitz conditions (see assumption \ref{standard assumptions}), let us again assume that $C$ is bounded. Furthermore we assume that all coefficients are continuously differentiable with respect to $t$. If $\tilde{C}\neq 0$ we also assume that $\bar{P}_0$ is invertible and that \eqref{regularity assumption} holds. Then there exists a unique solution $\bar{X}$ of \eqref{nonlinear mf ensemble kalman filter - true inverse} on the time interval $[0,T]$.
	\end{Thm}
	\begin{proof}
		
		As we have already pointed out before, for general $H$ and $R$, we can not expect the trace of the covariance matrix $\mathrm{tr}~\mathbb{Cov}_{Y}\left[X^P\right]$ to be bounded independent of $P$.\newline
		
		To make up for this, we partially localize the dynamics in a manner that only depends on the observations $Y$. To this end define for arbitrary $k\in\N$ 
		\begin{align}\label{definition psi}
			\tilde{\indicator}_k:\R\to[0,1]:x\mapsto\left(\indicator_{\left[-1/2,k-1/2]\right]}*\rho\right)(x),
		\end{align}
		where $\rho$ is a standard mollifier.\newline
		
		With this we define the process $\left(\bar{X}^{k}_t\right)_{t\in[0,T]}$ to be the solution of
		\begin{align}\label{localized equation}
			\begin{split}
				\mathrm{d}\bar{X}^k_t
				&=
				B_t(\bar{X}^k_t)\mathrm{d}t+
				C_t(\bar{X}^k_t)\mathrm{d}\bar{W}_t+
				\tilde{C}_t\mathrm{d}\bar{V}_t
				\\&\phantom{=}+
				\bar{\indicator}_k\left(\bar{P}^k_t H^\mathrm{T}_t+\tilde{C}_t\Gamma^{\mathrm{T}}_t\right) 
				R^{-1}_t\left(\mathrm{d}Y_t-\frac{H_t\left(\bar{X}^k_t+\bar{m}^k_t\right)}{2}\mathrm{d}t\right)
				\\&\phantom{=}-
				\bar{\indicator}_k\frac{\bar{P}^k_t H^\mathrm{T}_t+\tilde{C}_t\Gamma^{\mathrm{T}}_t}{2}
				R^{-1}_t \Gamma_t \tilde{C}^{\mathrm{T}}_t \tilde{\indicator}_{-k}\left(\bar{P}^k_t\right)^{-1} \left(\bar{X}^k_t-\bar{m}^k_t\right) \mathrm{d}t,
			\end{split}
		\end{align}	
		where
		\begin{align*}
			\bar{m}^k_t:=\mathbb{E}_{Y}
			\left[\bar{X}^k_t\right],~
			\bar{P}^k_t
			&:=\mathbb{Cov}_{Y}
			\left[\bar{X}^k_t\right],~
			\bar{\indicator}_{k}:=\tilde{\indicator}_k\left(|\bar{P}^k_t|^2\right)~\text{and}~
			\bar{\indicator}_{-k}
			:=\tilde{\indicator}_k\left(\left|\left(\bar{P}^k_t\right)^{-1}\right|^2\right).
		\end{align*}

		Note that \eqref{localized equation} still falls outside of the standard framework for the analysis of McKean--Vlasov equations (found e.g. in \cite{CarDel}), as the product of a bounded and an unbounded Lipschitz function may still not be Lipschitz, and thus the coefficients in \eqref{localized equation} therefore are only locally Lipschitz.\newline
		
		However the existence and uniqueness of such a solution $\bar{X}^k$ for every $k\in\N$ can be proven just as in the proof of Lemma \ref{Existence and uniqueness for mf - special case} by making a fixed point argument with respect to the covariance matrix $\bar{P}^k$, as the involvement of $\tilde{\indicator}_k$ bounds the covariance matrix independent of the argument in the fixed point map. Let $\bar{X}^{P,k}$ be the unique solution of equation
		\begin{align*}
			\begin{split}
				\mathrm{d}\bar{X}^{P,k}_t
				&=
				B_t\left(\bar{X}^{P,k}_t\right)\mathrm{d}t+
				C_t\left(\bar{X}^{P,k}_t\right)\mathrm{d}\bar{W}_t+
				\tilde{C}_t\mathrm{d}\bar{V}_t
				\\&\phantom{=}+
				\tilde{\indicator}_k(|P_t|^2) \left(P_t H^\mathrm{T}_t+\tilde{C}_t \Gamma^{\mathrm{T}}_t\right) 
				R^{-1}_t\left(\mathrm{d}Y_t-\frac{H_t\left(\bar{X}^{P,k}+\bar{m}^{P,k}_t\right)}{2}\mathrm{d}t\right)
				\\&\phantom{=}-
				\tilde{\indicator}_k(|P_t|^2)\left(P_t H^\mathrm{T}_t+\tilde{C}_t \Gamma^{\mathrm{T}}_t \right)R^{-1}_t
				\Gamma_t \tilde{C}^{\mathrm{T}}_t \tilde{\indicator}_k\left(\left|\left(P_t\right)^{-1}\right|^2\right)P_t^{-1} \frac{\bar{X}^{P,k}-\bar{m}^{P,k}_t}{2}\mathrm{d}t,
			\end{split}
		\end{align*}
		for any given $P\in\mathrm{SPSD}^1_Y([0,T])$.\newline
		
		In this case we have
		\begin{align}\label{cov XPk}
			\begin{split}
				\frac{\mathrm{d}\mathbb{Cov}_{Y}\left[\bar{X}^{P,k}_t\right]}{\mathrm{d}t}
				&=
				\comCov{B_t}{\bar{X}^{P,k}_t}
				-
				\tilde{\indicator}_k(|P_t|^2)\left(P_t H^\mathrm{T}_t+\tilde{C}_t\right)
				R^{-1}_t
				\frac{H_t\mathbb{Cov}_{Y}\left[\bar{X}^{P,k}_t\right]}{2}
				\\&\phantom{=}-
				\frac{\mathbb{Cov}_{Y}\left[\bar{X}^{P,k}_t\right] H^\mathrm{T}_t}{2}  R^{-1}_t
				\left(H_t P_t + \Gamma_t\tilde{C}^{\mathrm{T}}_t \right)\tilde{\indicator}_k(|P_t|^2)
				\\&\phantom{=}-
				\tilde{\indicator}_k(|P_t|^2)\frac{P_t H^\mathrm{T}_t+\tilde{C}_t\Gamma^{\mathrm{T}}_t}{2} R^{-1}_t \Gamma_t \tilde{C}^{\mathrm{T}}_t \tilde{\indicator}_k(|P^{-1}_t|^2)P^{-1}_t\mathbb{Cov}_{Y}\left[\bar{X}^{P,k}_t\right]
				\\&\phantom{=}-
				\mathbb{Cov}_{Y}\left[\bar{X}^{P,k}_t\right] P^{-1}_t \tilde{\indicator}_k(|P^{-1}_t|^2)
				\tilde{C}_t \Gamma^{\mathrm{T}}_t R^{-1}_t\frac{H_t\bar{P}_t +\Gamma_t\tilde{C}^{\mathrm{T}}_t}{2}\tilde{\indicator}_k(|P_t|^2)
				\\&\phantom{=}+
				\mathbb{E}_{Y}\left[C\left(\bar{X}^{P,k}_t\right)C\left(\bar{X}^{P,k}_t\right)^{\mathrm{T}}\right]
				+
				\tilde{C}\tilde{C}^{\mathrm{T}}
			\end{split}
		\end{align}
		
		Since  $\indicator_{\left[0,k-1]\right]}\leq\tilde{\indicator}_k\leq\indicator_{\left[-1,k]\right]}$, it holds that
		\begin{align*}
			\left|\tilde{\indicator}_k(|P_t|^2)\left(P_t H^\mathrm{T}_t+\tilde{C}_t \Gamma^{\mathrm{T}}_t\right).\right|
			&\leq
			\sqrt{k}\left|H^\mathrm{T}_t\right|+\left|\tilde{C}_t\Gamma^{\mathrm{T}}_t\right|
			\\
			\left|\tilde{\indicator}_k(|P^{-1}_t|^2)P^{-1}_t\right|
			&\leq
			\sqrt{k}.
		\end{align*}
		and one can derive for every fixed $k$ the boundedness of $\mathrm{tr}~\mathbb{Cov}_{Y}\left[\bar{X}^{P,k}_t\right]$ independent of $P$ by using the Grönwall inequality. As both
		\begin{align*}
			P\mapsto \tilde{\indicator}_k(|P|^2) \left(P H^\mathrm{T}+\tilde{C}\Gamma^{\mathrm{T}}_t\right) 
			~\text{and}~
			P\mapsto
			\tilde{\indicator}_k(|P|^2)\frac{P H^\mathrm{T}+\tilde{C}\Gamma^{\mathrm{T}}_t}{2} R^{-1}_t\Gamma_t \tilde{C}^{\mathrm{T}} \tilde{\indicator}_k\left(\left|\left(P\right)^{-1}\right|^2\right)P^{-1}
		\end{align*}
		are smooth functions with compact support (and therefore Lipschitz) one can now derive the existence and uniqueness of $\bar{X}^k$ for every fixed $k$ just as in the proof of Lemma \ref{Existence and uniqueness for mf - special case}.\newline

		For the fixed point $\bar{X}^k$ equation \eqref{cov XPk} gives us
		\begin{align}\label{cov Xk}
			\begin{split}
				\frac{\mathrm{d}\mathrm{tr} \bar{P}^k_t}{\mathrm{d}t}
				&=
				\comCov{B_t}{\bar{X}^k_t}
				-
				\bar{\indicator}_{k}
				\bar{P}^k_t H^{\mathrm{T}}_t R^{-1}_t H \bar{P}^k_t
				-
				\bar{\indicator}_{k} \left(1+\bar{\indicator}_{-k}\right)
				\mathrm{Sym}\left(\tilde{C}_t \Gamma^{\mathrm{T}}_t R^{-1}_t H_t \bar{P}^k_t\right)
				\\&\phantom{=}+
				\mathbb{E}_{Y}\left[C\left(X^{P,k}_t\right)C\left(X^{P,k}_t\right)^{\mathrm{T}}\right]
				+
				\tilde{C}_t
				\left(
				1-
				\bar{\indicator}_{k}
				\bar{\indicator}_{-k} \Gamma^{\mathrm{T}}_t R^{-1}_t \Gamma_t
				\right) \tilde{C}^{\mathrm{T}}_t
			\end{split}
		\end{align}
		and thus one derives just as in subsection \ref{covariance bounds}
		\begin{align*}
			\begin{split}
				\frac{\mathrm{d}\mathrm{tr} \bar{P}^k_t}{\mathrm{d}t}
				&\leq
				\left(2\mathrm{Lip}\left(B\right)+d_x \bar{\indicator}_{k}\left(1+\bar{\indicator}_{-k}\right) \left|\tilde{C}_t \Gamma^{\mathrm{T}}_t R^{-1}_t H_t\right|\right)
				~\mathrm{tr} \bar{P}^k_t
				\\&\phantom{=}+
				\norm{C_t}{\infty}^2
				+
				\left|\tilde{C}_t\right|^2
				+ 
				d_x
				\bar{\indicator}_{k}
				\bar{\indicator}_{-k}
				\left|\tilde{C}_t \Gamma^{\mathrm{T}}_t R^{-1}_t \Gamma_t\tilde{C}^{\mathrm{T}}_t\right|
				\\&\leq
				2\left(\mathrm{Lip}\left(B\right)+d_x \left|\tilde{C}_t \Gamma^{\mathrm{T}}_t R^{-1}_t H_t\right|\right)
				~\mathrm{tr} \bar{P}^k_t
				+
				\norm{C_t}{\infty}^2
				+
				\left|\tilde{C}_t\right|^2
				+ 
				d_x\left|\tilde{C}_t \Gamma^{\mathrm{T}}_t R^{-1}_t \Gamma_t\tilde{C}^{\mathrm{T}}_t\right|,
			\end{split}
		\end{align*}
		which is similar to the differential inequality \eqref{trace inequality - dynamical}. Therefore we derive by the Grönwall lemma that 
		\begin{align*}
			\sup_{t\leq T}\mathrm{tr}\bar{P}^k_t
			&\leq \exp(2\int_{0}^{T}d_x \mathrm{Lip}(B)+\left|\tilde{C}_t \Gamma^{\mathrm{T}}_t R^{-1}_t H_t\right|\mathrm{d}t)
			\\&\phantom{=}\times\left(\bar{P}_0+\int_{0}^{T}\norm{C_t}{\infty}^2
			+
			\left|\tilde{C}_t\right|^2+d_x\left|\tilde{C}_t \Gamma^{\mathrm{T}}_t R^{-1}_t \Gamma_t\tilde{C}^{\mathrm{T}}_t\right|\mathrm{d}t\right)
			:=\hat{\Psi}(T)
		\end{align*}
		(the constant is defined in \eqref{a priori bound covariance}) and we see that if 
		\begin{align*}
			k
			\geq
			\hat{\Psi}(T)^2+1,
		\end{align*}
		then $\bar{\indicator}_{k}=\tilde{\indicator}_k(|\bar{P}^k_t|^2)=1$ for all $t\in[0,T]$.\newline

		Similarly we can bound $\lambda_{\min}\left(\bar{P}^k_t\right)$ from below. Using \eqref{cov Xk}, we see, by employing the same inequalities as in Lemma \ref{regularity of covariance matrix}, that for the $i$-th eigenvalue $\lambda^{k,i}$ of $\bar{P}^k$, the differential inequality 
		\begin{align*}
			\begin{split}
				\frac{\mathrm{d}\lambda^{k,i}_t}{\mathrm{d}t}
				&\geq
				-2\mathrm{Lip}(B)~\sqrt{\bar{\Psi}(T)} \sqrt{\lambda^{k,i}_t}
				-
				\bar{\indicator}_{k}
				\lambda_{\max}\left(H^\mathrm{T}_t R^{-1}_t H_t\right) \left(\lambda^{i}_t\right)^2
				\\&\phantom{=} 
				-2 \bar{\indicator}_{k} \left(1+\bar{\indicator}_{-k}\right)
				\lambda_{\max}\left(\mathrm{Sym}\left(\tilde{C}_t \Gamma^{\mathrm{T}}_t R^{-1}_t H_t\right)\right)
				\lambda^i_t
				\\&\phantom{=}+
				\inf_{x}\lambda_{\min}\left(C_t(x)C_t(x)^{\mathrm{T}}\right)
				+
				\lambda_{\min}\left(\tilde{C}_t \left(I-\bar{\indicator}_{k}\bar{\indicator}_{-k}\Gamma^{\mathrm{T}}_t R^{-1}_t \Gamma_t\right) \tilde{C}^{\mathrm{T}}_t\right)
			\end{split}
		\end{align*}
		holds. Since $0\leq \bar{\indicator}_{k},\bar{\indicator}_{-k}\leq 1$, we have, due to the variational characterization of the smallest eigenvalue,
		\begin{align*}
			&\lambda_{\min}\left(\tilde{C}_t \left(I-\bar{\indicator}_{k}\bar{\indicator}_{-k} \Gamma^{\mathrm{T}}_t R^{-1}_t \Gamma_t\right) \tilde{C}^{\mathrm{T}}_t\right)
			=
			\min_{|v|=1} 
			v^{\mathrm{T}}\tilde{C}_t\tilde{C}^{\mathrm{T}}_t v
			-
			\bar{\indicator}_{k}\bar{\indicator}_{-k}v^{\mathrm{T}}\tilde{C}_t \Gamma^{\mathrm{T}}_t R^{-1}_t \Gamma_t\tilde{C}^{\mathrm{T}}_t v
			\\&\phantom{=}\geq
			\min_{|v|=1} 
			v^{\mathrm{T}}\tilde{C}_t\tilde{C}^{\mathrm{T}}_t v
			-
			v^{\mathrm{T}}\tilde{C}_t \Gamma^{\mathrm{T}}_t  R^{-1}_t \Gamma_t\tilde{C}^{\mathrm{T}}_t v
			\geq
			\lambda_{\min}\left(\tilde{C}_t \left(I-\Gamma^{\mathrm{T}}_t R^{-1}_t \Gamma_t\right) \tilde{C}^{\mathrm{T}}_t\right).
		\end{align*}
		
		This lets us derive the same inequality as in the proof of Lemma \ref{regularity of covariance matrix}
		\begin{align*}
			\begin{split}
				\frac{\mathrm{d}\lambda^{k,i}_t}{\mathrm{d}t}
				&\geq
				-2\mathrm{Lip}(B)~\sqrt{\bar{\Psi}(T)} \sqrt{\lambda^{k,i}_t}
				\\&\phantom{=}-
				\lambda_{\max}\left(H^\mathrm{T}_t R^{-1}_t H_t\right) \left(\lambda^{i}_t\right)^2
				-4
				\lambda_{\max}\left(\mathrm{Sym}\left(\tilde{C}_t\Gamma^{\mathrm{T}}_t R^{-1}_t H_t\right)\right)
				\lambda^i_t
				+\gamma_t.
			\end{split}
		\end{align*}
		
		This inequality allows us to bound $\lambda_{\min}\left(\bar{P}^k_t\right)$ from below, independently of $k$. To see this, let $\underline{\lambda}^i$ be a solution of the initial value problem
		\begin{align*}
			\begin{cases}
				\frac{\mathrm{d}\underline{\lambda}^i_t}{\mathrm{d}t}
				&=
				-2\mathrm{Lip}(B)~\sqrt{\bar{\Psi}(T)} \sqrt{\underline{\lambda}^i_t}
				\\&\phantom{=}-
				\lambda_{\max}\left(H^\mathrm{T}_t R^{-1}_t H_t\right) \left(\underline{\lambda}^i_t\right)^2
				-4
				\lambda_{\max}\left(\mathrm{Sym}\left(\tilde{C}_t\Gamma^{\mathrm{T}}_t R^{-1}_t H_t\right)\right)
				\underline{\lambda}^i_t
				+
				\frac{\gamma_t}{2}
				\\
				\underline{\lambda}^i_0&=\frac{\bar{\lambda}^i_0}{2},
			\end{cases}
		\end{align*}
		where $\bar{\lambda}^i_0$ denotes the $i$-th eigenvector of $\bar{P}_0$. Note that the Peano theorem guarantees the existence of such a solution, while its uniqueness is not guaranteed. Thus we have $\left(\lambda^{k,i}\right)'(t_0)>\left(\underline{\lambda}^i\right)'(t_0)$, if $\lambda^{k,i}_{t_0}=\bar{\lambda}^i_{t_0}$. Now \cite[II. Lemma, page 64]{Walter} guarantees that $\lambda^{k,i}_t\geq\bar{\lambda}^i_t$ and furthermore it also guarantees that $\bar{\lambda}^i_t>0$ for all times $t$. Therefore we have found a uniform (both in $k$ and in $i$) lower bound.\newline
		
		This lets us conclude that for $k\in\N$ large enough such that $k\geq\frac{1}{d_x \left(\min_{t\leq T}\bar{\lambda}^i_t\right)^2}$, we have $\indicator_{k}\left(\left|\left(\bar{P}^k_t\right)^{-1}\right|^2\right)=1$ and thus $\bar{X}^k$ is the unique solution of \eqref{nonlinear mf ensemble kalman filter - true inverse} on $[0,T]$.\newline

	\end{proof}

	\subsection{Generalization to other EnKBFs}\label{section other EnKBF}
	
	Note that for the uncorrelated case $\tilde{C}=0$ the solution of \eqref{nonlinear mf ensemble kalman filter - true inverse} is not the only McKean--Vlasov equation that describes the mean field limit of an EnKBF. Indeed in the literature it is sometimes referred to as the deterministic (mean field) EnKBF \cite{Riccati diff - stability}\cite{CriDMJasRuz}. This is the continuous time counterpart of the filter defined by Sakov and Oke in \cite{Sakov Oke} and was studied  in \cite{Riccati diff - stability}, \cite{Riccati diff - Perturbation} as well as in \cite{Reich Cotter}. Given a fixed realization/path of the observations $Y$, it is a deterministic perturbation of the signal. However, even for a fixed path of $Y$, it does not define a deterministic equation due to the involvement of the Brownian motions $\bar{V},\bar{W}$. In filtering one often replaces Brownian motions with deterministic inflation terms \cite{Reich Cotter}\cite{DeWiljesReichStannat} to derive deterministic equations. This is based on the following observation that was derived in \cite{Reich Schroedinger}.
	\begin{Lemma}\label{Lemma log SDE}
		Let $Z$ be the solution of the SDE
		\begin{align}\label{general SDE}
			\mathrm{d}Z_t=\alpha_t(Z_t)\mathrm{d}t+\beta_t\mathrm{d}W_t.
		\end{align}
		For every time $t$ denote by $\eta^Z_t$ the density of the law of $Z_t$. Furthermore we assume that there exists a solution $\tilde{Z}$ to
		\begin{align*}
			\mathrm{d}\tilde{Z}_t=\alpha_t(\tilde{Z}_t)\mathrm{d}t-\frac{\beta_t\beta_t^{\mathrm{T}}}{2}\nabla\log\eta^{\tilde{Z}}_t(\tilde{Z})\mathrm{d}t,
		\end{align*}
		where $\eta^{\tilde{Z}}_t$ denotes the density of the law of $\tilde{Z}_t$. We assume that both $\eta^Z,\eta^{\tilde{Z}}$ are well defined densities and that both are strictly positive and smooth.  Then if the Kolmogorov backward equation 
		\begin{align}\label{general Kolmogorov backward}
			\partial_t\eta^Z_t=\mathrm{div}\left(\frac{\beta_t\beta_t^{\mathrm{T}}}{2}\nabla\eta^Z_t\right)-\mathrm{div}\left(\alpha_t\eta^Z_t\right)
		\end{align}
		associated to \eqref{general SDE} has a unique (strong) solution, it holds that $\eta^Z_t=\eta^{\tilde{Z}}_t$.
	\end{Lemma}
	\begin{proof}
		Assuming that $\eta^{\tilde{Z}}$ is sufficiently smooth, it satisfies
		\begin{align*}
			\partial_t\eta^{\tilde{Z}}_t
			=
			-\mathrm{div}\left(\left(-\frac{\beta_t\beta_t^{\mathrm{T}}}{2}\nabla\log\eta^{\tilde{Z}}_t+\alpha_t\right)\eta^{\tilde{Z}}_t\right)
			=
			\mathrm{div}\left(\frac{\beta_t\beta_t^{\mathrm{T}}}{2}\nabla\eta^{\tilde{Z}}_t\right)-\mathrm{div}\left(\alpha_t\eta^{\tilde{Z}}_t\right),
		\end{align*}
		and thus by the uniqueness of solutions to \eqref{general Kolmogorov backward}, we derive $\eta^{Z}_t\eta^{\tilde{Z}}_t$ for all $t\geq 0$.
	\end{proof}

	If we now assume that for every point in time $t$ the matrix $C_t$ does not depend on state $\bar{X}_t$, we could thus replace $C_t\mathrm{d}\bar{W}_t + \tilde{C}_t\mathrm{d}\bar{V}_t$ in \eqref{nonlinear mf ensemble kalman filter - true inverse} by $-\frac{C_tC_t^{\mathrm{T}}+\tilde{C}_t\tilde{C}_t^{\mathrm{T}}}{2}\nabla\log\bar{\eta}_t$. If one uses \eqref{projection relation} to approximate $\nabla\log\bar{\eta}_t$ by its affine projection $\pi^1[\bar{\eta}_t]\nabla\log\bar{\eta}_t=-\bar{P}^{-1}_t\left(\cdot-\bar{m}_t\right)$, one thus derives the equation
	\begin{align}\label{mf limit of transport EnKBF}
		\begin{split}
			\mathrm{d}\bar{X}_t
			&=
			B_t(\bar{X}_t)\mathrm{d}t
			+
			\left(C_t C^{\mathrm{T}}_t+\tilde{C}_t \tilde{C}^{\mathrm{T}}_t\right) \bar{P}_t^{-1}\frac{\bar{X}_t-\bar{m}_t}{2}\mathrm{d}t
			\\&\phantom{=}+
			\left(\bar{P}_t H^{\mathrm{T}}_t+\tilde{C}_t\Gamma^{\mathrm{T}}_t\right)R^{-1}_{t}
			\left(\mathrm{d}Y_t-\frac{H_t\left(\bar{X}_t+\bar{m}_t\right)}{2}\mathrm{d}t\right)
			\\&\phantom{=}-
			\left(\bar{P}_t H^{\mathrm{T}}+\tilde{C}_t\Gamma^{\mathrm{T}}_t\right)R^{-1}_{t}\Gamma_t
			\tilde{C}^{\mathrm{T}}_t\bar{P}^{-1}_t\frac{\bar{X}_t-\bar{m}_t}{2}\mathrm{d}t.
		\end{split}
	\end{align}

	This correspondence of Brownian motion and the inflation term $\bar{P}^{-1}\left(\cdot-\bar{m}_t\right)$ can also be used in the other direction to derive a generalization of the classical EnKBF with randomized innovation term, that also eliminates the singular inflation term which stems from the correlated observation noise. To this end we note that \eqref{mf limit of transport EnKBF} can be rewritten as
	\begin{align*}
		\begin{split}
			\mathrm{d}\bar{X}_t
			&=
			B_t(\bar{X}_t)\mathrm{d}t
			+
			\left(C_t C^{\mathrm{T}}_t+\tilde{C}_t \tilde{C}^{\mathrm{T}}_t\right) \bar{P}_t^{-1}\frac{\bar{X}_t-\bar{m}_t}{2}\mathrm{d}t
			\\&\phantom{=}+
			\left(\bar{P}_t H^{\mathrm{T}}_t+\tilde{C}_t\Gamma^{\mathrm{T}}_t\right)R^{-1}_{t}
			\left(\mathrm{d}Y_t-H_t\bar{X}_t\mathrm{d}t\right)
			\\&\phantom{=}+
			\left(\bar{P}_t H^{\mathrm{T}}_t+\tilde{C}_t\Gamma^{\mathrm{T}}_t\right)R^{-1}_{t}
			\left(H_t\bar{P}_t-\Gamma_t\tilde{C}^{\mathrm{T}}_t\right)\bar{P}^{-1}_t\frac{\bar{X}_t-\bar{m}_t}{2}\mathrm{d}t.
		\end{split}
	\end{align*}
	
	Again we use \eqref{projection relation} to replace $-\bar{P}^{-1}_t\left(\bar{X}_t-\bar{m}_t\right)$ by $\nabla\log\bar{\eta}\left(\bar{X}_t\right)$ and thus turn \eqref{mf limit of transport EnKBF} into
	\begin{align}\label{mf equation all log}
		\begin{split}
			\mathrm{d}\bar{X}_t
			&=
			B_t(\bar{X}_t)\mathrm{d}t
			-
			\frac{C_t C^{\mathrm{T}}_t+\tilde{C}_t \tilde{C}^{\mathrm{T}}_t}{2} \nabla\log\bar{\eta}\left(\bar{X}_t\right)\mathrm{d}t
			\\&\phantom{=}+
			\left(\bar{P}_t H^{\mathrm{T}}_t+\tilde{C}_t\Gamma^{\mathrm{T}}_t\right)R^{-1}_{t}
			\left(\mathrm{d}Y_t-H_t\bar{X}_t\mathrm{d}t\right)
			\\&\phantom{=}-
			\frac{\left(\bar{P}_t H^{\mathrm{T}}_t+\tilde{C}_t\Gamma^{\mathrm{T}}_t\right)R^{-1}_{t}
				\left(H_t\bar{P}_t-\Gamma_t\tilde{C}^{\mathrm{T}}_t\right)}{2}
			\nabla\log\bar{\eta}\left(\bar{X}_t\right)\mathrm{d}t.
		\end{split}
	\end{align}
	
	Next we note that since the matrix $\mathcal{M}:=\tilde{C}_t\Gamma^{\mathrm{T}}_tR^{-1}_{t}H_t\bar{P}_t-\bar{P}_t H^{\mathrm{T}}_tR^{-1}_{t}\Gamma_t\tilde{C}^{\mathrm{T}}_t$ is skew symmetric, we have
	\begin{align*}
		\mathrm{div}\left(\bar{\eta}_t \mathcal{M}\nabla\log\bar{\eta}\right)
		&=
		\mathrm{div}\left( \mathcal{M}\nabla\bar{\eta}\right)
		=
		\mathcal{M}:\bar{\eta}_t''
		=\mathrm{tr}\left[\mathcal{M}^{\mathrm{T}}\bar{\eta}_t''\right]
		=\mathrm{tr}\left[\bar{\eta}_t''\mathcal{M}\right]
		=\mathrm{tr}\left[\mathcal{M}\bar{\eta}_t''\right]
		\\&=
		-\mathrm{tr}\left[\mathcal{M}^{\mathrm{T}}\bar{\eta}_t''\right]
		=
		-\mathrm{div}\left(\bar{\eta}_t \mathcal{M}\nabla\log\bar{\eta}\right)
		=
		0.
	\end{align*}
	
	Thus we can erase $\mathcal{M}=\tilde{C}_t\Gamma^{\mathrm{T}}_tR^{-1}_{t}H_t\bar{P}_t-\bar{P}_t H^{\mathrm{T}}_tR^{-1}_{t}\Gamma_t\tilde{C}^{\mathrm{T}}_t$ in \eqref{mf equation all log} without changing the associated Kolmogorov equation. Therefore solutions of
	\begin{align}\label{mean field vanilla with logs}
		\begin{split}
			\mathrm{d}\bar{X}_t
			&=
			B_t(\bar{X}_t)\mathrm{d}t
			-
			\frac{C_t C^{\mathrm{T}}_t}{2} \nabla\log\bar{\eta}\left(\bar{X}_t\right)\mathrm{d}t
			\\&\phantom{=}+
			\left(\bar{P}_t H^{\mathrm{T}}_t+\tilde{C}_t\Gamma^{\mathrm{T}}_t\right)R^{-1}_{t}
			\left(\mathrm{d}Y_t-H_t\bar{X}_t\mathrm{d}t\right)
			\\&\phantom{=}-
			\frac{\bar{P}_t H^{\mathrm{T}}_tR^{-1}_{t}H_t\bar{P}_t
				+
				\tilde{C}_t\left(I-\Gamma^{\mathrm{T}}_t R^{-1}_{t}\Gamma_t\right)\tilde{C}^{\mathrm{T}}_t
			}{2}
			\nabla\log\bar{\eta}\left(\bar{X}_t\right)\mathrm{d}t.
		\end{split}
	\end{align}
	have the same conditional laws $\left(\bar{\eta}_t\right)_{t\geq 0}$ as solutions to \eqref{mf equation all log}. Since $\Gamma^{\mathrm{T}}_t R^{-1}_{t}\Gamma_t$ is a projection onto the range of $\Gamma^{\mathrm{T}}_t$, it is easy to see that we can rewrite \eqref{mean field vanilla with logs} into
	\begin{align*}
		\begin{split}
			&\mathrm{d}\bar{X}_t
			=
			B_t(\bar{X}_t)\mathrm{d}t
			-
			\frac{C_t C^{\mathrm{T}}_t}{2} \nabla\log\bar{\eta}\left(\bar{X}_t\right)\mathrm{d}t
			+
			\left(\bar{P}_t H^{\mathrm{T}}_t+\tilde{C}_t\Gamma^{\mathrm{T}}_t\right)R^{-1}_{t}
			\left(\mathrm{d}Y_t-H_t\bar{X}_t\mathrm{d}t\right)
			\\&-
			\frac{\left(\tilde{C}_t\left(I-\Gamma^{\mathrm{T}}_tR^{-1}_t\Gamma_t\right)-\bar{P}_t H^{\mathrm{T}}_tR^{-1}_{t}\Gamma_t\right)
				\left(\left(I-\Gamma^{\mathrm{T}}_tR^{-1}_t\Gamma_t\right)\tilde{C}^{\mathrm{T}}_t- 
				\Gamma^{\mathrm{T}}_t R^{-1}_t	H_t	\bar{P}_t \right)
			}{2}
			\nabla\log\bar{\eta}\left(\bar{X}_t\right)\mathrm{d}t.
		\end{split}
	\end{align*}
	
	By now applying Lemma \ref{Lemma log SDE} we derive the equation
	\begin{align}\label{vanilla EnKBF}
		\begin{split}
			\mathrm{d}\bar{X}_t
			&=
			B_t(\bar{X}_t)\mathrm{d}t
			+
			C_t \mathrm{d}\bar{W}_t
			+
			\tilde{C}_t \mathrm{d}\bar{V}_t
			\\&\phantom{=}+
			\left(\bar{P}_t H^{\mathrm{T}}_t+\tilde{C}_t\Gamma^{\mathrm{T}}_t\right)
			R^{-1}_{t}
			\left(\mathrm{d}Y_t-H_t\bar{X}_t\mathrm{d}t -\Gamma_t\mathrm{d}\bar{V}_t\right).
		\end{split}
	\end{align}
	{
		This is a generalization of the classical EnKBF with randomized innovation term to the correlated noise case. This mean field equation was also derived by different arguments in \cite{NueReiRoz} for state and parameter estimation problems. Note that \eqref{vanilla EnKBF} does not contain any singular terms. This simplifies the analysis and in \cite{Coghi et al} a well posedness proof for \eqref{vanilla EnKBF} can be found under the restriction that the observation function $H$ must be bounded. Our well posedness result of \eqref{nonlinear mf ensemble kalman filter - true inverse}  in Theorem \ref{Existence and uniqueness for mf} can easily be adapted to all McKean--Vlasov equations in this section, in particular to the mean field filters \eqref{vanilla EnKBF} and \eqref{mf limit of transport EnKBF}.}
	
	\begin{Rmk}
		We note that since for Gaussian $\bar{\eta}_t$ it holds that $\pi^1[\bar{\eta}_t]\nabla\log\bar{\eta}_t= \nabla\log\bar{\eta}_t$, and by Lemma \ref{Lemma log SDE}, all mean field filters derived in this subsection are optimal in the linear Gaussian setting. In particular \eqref{vanilla EnKBF} and \eqref{mf limit of transport EnKBF} are consistent representations of the true posterior in this case. 
	\end{Rmk}
	
	\section{The well posedness of the EnKBF}\label{section particle system}

	Just as for the linear Gaussian case, one can approximate \eqref{nonlinear mf ensemble kalman filter - true inverse} in a straight forward manner by the interacting particle system 
	\begin{align}\label{nonlinear Ensemble Kalman}
		\begin{split}
			\mathrm{d}X^i_t
			&=
			B_t\left(X^i_t\right)\mathrm{d}t+
			C_t\left(X^i_t\right)\mathrm{d}W^i_t
			+
			\tilde{C}_t\mathrm{d}V^i_t
			\\&\phantom{=}+
			\left(P^M_t H^\mathrm{T}_t+\tilde{C}_t\Gamma^{\mathrm{T}}_t\right)
			R^{-1}_t
			\left(\mathrm{d}Y_t-\frac{H_t\left(X^i_t+x^M_t\right)}{2}\mathrm{d}t\right)
			\\&\phantom{=}-
			\left(	P^M_t H^\mathrm{T}_t+\tilde{C}_t\Gamma^{\mathrm{T}}_t\right)~
			R^{-1}_t~\Gamma_t
			\tilde{C}^{\mathrm{T}}_t\left(P^M_t\right)^{+}
			\frac{X^i_t-x^M_t}{2}\mathrm{d}t
		\end{split}
	\end{align}
	for $i=1,\cdots,M$. Where again 
	\begin{align*}
		x^M_t:=\frac{1}{M}\sum_{i=1}^{M} X^i_t
		~\text{and}~
		P^M_t:=\frac{1}{M-1}\sum_{i=1}^{M}\left(X^i_t-x^M_t\right)\left(X^i_t-x^M_t\right)^\mathrm{T}
	\end{align*}
	denote the ensemble average and the ensemble covariance matrix.  $\left(P^M_t\right)^{+}$ is the Moore--Penrose pseudoinverse of $P^M_t$.\newline
	
	In the uncorrelated case $\tilde{C}=0$ system \eqref{nonlinear Ensemble Kalman} is referred to as the (deterministic) EnKBF and we will also refer to it as such in the correlated noise framework. Only in the linear Gaussian setting it will be an approximation of the optimal filter, otherwise its mean field limit \eqref{nonlinear mf ensemble kalman filter - true inverse} will not be a representation of the posterior.\newline

	Equation \eqref{nonlinear Ensemble Kalman} is a system of nonlinear SDEs, in which the coefficients do not satisfy linear growth properties and which therefore falls outside the standard existence theory for SDEs found for example in \cite{KaratzasShreve}. Nevertheless in \cite{StannatLange} the existence and uniqueness of solutions was proven for uncorrelated observation noise. More precisely it was assumed that $C$ is constant and $\tilde{C}=0$.\newline
	
	However, generalizing this result to the correlated case \eqref{nonlinear Ensemble Kalman} is not straightforward, as the pseudoinverse $\left(P^M_t\right)^{+}$ does not depend continuously on the ensemble members $X^i,~i=1,\cdots,M$ and actually develops singularities where $P^M_t$ changes its rank. Indeed for $M=2$ it is easy to see that
	\begin{align}\label{inflation term 2 particles}
		\left(P^M_t\right)^{+}\left(X^i_t-x^M_t\right)
		=\indicator_{\left[X^1_t\neq X^2_t\right]}
		\frac{\left(X^i_t-x^M_t\right)}{2 \left|X^i_t-x^M_t\right|^2},~i=1,2,
	\end{align}
	which becomes singular when the two particles collide. This makes the analysis of \eqref{nonlinear Ensemble Kalman} particularly challenging as some of the most general existing well posedness results require at least local integrability. Indeed there exist counter examples of SDEs containing drift terms that are structurally similar to \eqref{inflation term 2 particles} that do not even admit weak solutions \cite[Example 1.17]{ChernyEngelbert}. The following theorem shows that despite these difficulties \eqref{nonlinear Ensemble Kalman} is well posed under suitable assumptions.

	\begin{Thm}\label{well posedness EnKBF theorem}
		Beside the standard Lipschitz assumptions \ref{standard assumptions}, let us assume that $C$. If $\tilde{C}\neq 0$ we furthermore assume that $P^M_0$ is invertible and that $M$ is large enough, so that
		\begin{align}\label{condtion gammaM}
			\underline{\gamma}^M:=
			\inf_{t\leq T}
			&
			\left(
			\gamma_t-
			\frac{2}{M-1}
			\left(1+\sqrt{d_x}\right) \left(\norm{C_t}{\infty}^2+\left|\tilde{C}_t\right|^2\right)
			\right)
			>
			0.
		\end{align}
		holds.	Then there exists a unique strong solution to \eqref{nonlinear Ensemble Kalman}.
	\end{Thm}
	
	\begin{Rmk}\label{M > dx}
		Before we prove Theorem \ref{well posedness EnKBF theorem} we want to remark that the assumption that $P^M_0$ is invertible implies $M\geq d_x+1$. Thus this result is restricted to large ensemble sizes.
	\end{Rmk}
	
	\begin{Rmk}\label{consistency of conditions}
		Note that in the limit $M\to\infty$ assumption \eqref{condtion gammaM} is consistent with condition \eqref{regularity assumption}, which was required for the well posedness of the mean field EnKBF \eqref{nonlinear mf ensemble kalman filter - true inverse} if $\tilde{C}\neq 0$.
	\end{Rmk}
	
	\begin{proof}[Proof of Theorem \ref{well posedness EnKBF theorem}]
		The proof consists of three major steps. First we derive equations for the empirical moments of the ensemble. Next we show that $P^M$ stays regular as long as it does not blow up. Finally we prove the boundedness of $P^M$ and by that also the existence and uniqueness of solutions.\newline
		
		Denote by $\bar{\xi}$ the explosion time of the ensemble and by $\underline{\xi}$ the time $P^M$ becomes non invertible. \newline

		\textbf{Step 1:} First we note that by the linearity of stochastic differentials the following equation for the mean holds up to time $\bar{\xi}\wedge\underline{\xi}$
		\begin{align*}
			\mathrm{d}x^M_t
			&=
			\underbrace{\frac{1}{M}\sum_{i=1}^{M}B_t(X^i_t)}_{:=b^M_t}\mathrm{d}t+
			\frac{1}{M}\sum_{i=1}^{M}C_t(X^i_t)\mathrm{d}W^i_t+
			\tilde{C}_t\mathrm{d}V^i_t
			\\&\phantom{=}+
			\left(P^M_t H^\mathrm{T}_t+\tilde{C}_t\Gamma^{\mathrm{T}}_t\right)
			R^{-1}_t
			\left(\mathrm{d}Y_t-H_t x^M_t\mathrm{d}t\right)
			,
		\end{align*}
		which gives us
		\begin{align*}
			\mathrm{d}\left(X^i_t-x^M_t\right)
			&=
			\left(B_t(X^i_t)-b^M_t\right)\mathrm{d}t
			\\&\phantom{=}+
			\frac{1}{M}\sum_{j=1}^{M} \left(C_t(X^i)\mathrm{d}W^i_t-C_t(X^j)\mathrm{d}W^j_t\right)
			+
			\frac{1}{M}\sum_{j=1}^{M} \left(\tilde{C}\mathrm{d}V^i_t-\tilde{C}\mathrm{d}V^j_t\right)
			\\&\phantom{=}-
			\left(P^M_t H^\mathrm{T}_t+\tilde{C}_t\Gamma^{\mathrm{T}}_t\right) R^{-1}_t \left(H_t+\Gamma_t\tilde{C}^{\mathrm{T}}_t\left(P^M_t\right)^{+}\right)\frac{X^i_t-x^M_t}{2}\mathrm{d}t.
		\end{align*}
		
		By employing Itô's product rule we derive the following evolution equation for the empirical covariance matrix up to time $\bar{\xi}\wedge\underline{\xi}$
		\begin{align}\label{evolution empirical covariances}
			\begin{split}
				\mathrm{d} P^{M}_t
				&=
				\comCov{B_t}{X_t}_M \mathrm{d}t
				+
				\left(\frac{1}{M}\sum_{i=1}^{M}C_t(X^i_t)C_t(X^i_t)^{\mathrm{T}}
				+
				\tilde{C}_t \tilde{C}_t^{\mathrm{T}}\right)
				\mathrm{d}t
				+\mathrm{d}\mathfrak{lm}_t
				\\&\phantom{=}-
				\mathrm{Sym}\left(
				P^M_t
				\left(H^{\mathrm{T}}_t+\left(P^M_t\right)^{+}\tilde{C}_t\Gamma^{\mathrm{T}}_t\right) R^{-1}_t
				\left(H_t P^M_t +\Gamma_t\tilde{C}^\mathrm{T}_t\right)\right)
				\mathrm{d}t,
			\end{split}
		\end{align}
		where
		\begin{align*}
			\comCov{B_t}{X_t}_M
			&:=
			\frac{2}{M-1}\sum_{i=1}^{M}
			\mathrm{Sym}
			\left(
			\left(B_t(X^i_t)-b^M_t\right)\left(X^i_t-x^M_t\right)^{\mathrm{T}}
			\right)
			\mathrm{d}t
			.
		\end{align*}
		and $\mathfrak{lm}$ denotes the local martingale starting at zero given by
		\begin{align}\label{local martingale emp cov}
			\begin{split}
				\mathrm{d}\mathfrak{lm}_t
				&:=
				\frac{2}{M-1}\sum_{j=1}^{M}
				\mathrm{Sym}
				\left(
				\left(X^j_t-x^M_t\right)
				\left(C(X^j_t)\mathrm{d}W^j_t\right)^{\mathrm{T}}
				\right)
				\\&\phantom{=}+
				\frac{2}{M-1}\sum_{j=1}^{M}
				\mathrm{Sym}
				\left(
				\left(X^j_t-x^M_t\right)
				\left(\tilde{C}\mathrm{d}V^j_t\right)^{\mathrm{T}}
				\right).
			\end{split}
		\end{align}

		\textbf{Step 2:} Using this equation we want to prove that $P^M$ is always invertible, as long as it stays bounded. To this end we define for every symmetric positive semidefinite matrix $P$ its regularized inverse $P^{+\epsilon}:=\left(P+\epsilon I\right)^{-1}$ for any $\epsilon>0$ and where $I$ denotes the identity matrix. Note that if $P$ is not invertible, then $\mathrm{tr}P^{+\epsilon}$ will blow up for $\epsilon\to 0$. Therefore we now try to dominate $\left(P^M\right)^{+\epsilon}$ uniformly in $\epsilon$.\newline
		
		In the following we will need the first and second order Fréchet derivatives of  $T_\epsilon(P):=P^{+\epsilon}$ which are given by
		\begin{align*}
			T_\epsilon'(P)[\mathcal{A}]&=-P^{+\epsilon} \mathcal{A} P^{+\epsilon}\\
			T_\epsilon''(P)[\mathcal{A},\mathcal{B}]&=
			2~\mathrm{Sym}\left(P^{+\epsilon} \mathcal{A} P^{+\epsilon} \mathcal{B} P^{+\epsilon}\right)
		\end{align*}
		for any two matrices $\mathcal{A}$ and $\mathcal{B}$. To simplify notations in the following we also define for any vector $\hat{x}$
		\begin{align}\label{definition mathcalS}
			\begin{split}
				\mathcal{S}(\mathcal{A},\hat{x},\mathcal{B})
				:=
				\sum_{k}
				\mathrm{tr}
				\left[\mathcal{A}
				\mathrm{Sym}\left(\hat{x} \mathcal{B}_{k}^{\mathrm{T}}\right)
				\mathcal{A}
				\mathrm{Sym}\left(\hat{x}\mathcal{B}_{k}^{\mathrm{T}}\right)
				\mathcal{A}\right],
			\end{split}
		\end{align}
		where $\mathcal{B}_{k}$ denotes the $k$-th column of $\mathcal{B}$.\newline
		
		Now we note that if we define the local martingale
		\begin{align*}
			\mathrm{d}\mathcal{M}^{\epsilon}_t
			&:=
			-\mathrm{tr}\left[\left(P^M_t\right)^{+\epsilon}~\mathrm{d}\mathfrak{lm}_t~\left(P^M_t\right)^{+\epsilon}\right],
		\end{align*}
		then by Itôs rule one can see
		\begin{align}
			\begin{split}
				\mathrm{d}\mathrm{tr}\left(P^M_t\right)^{+\epsilon}
				&=
				-\mathrm{tr}\left[\left(P^M_t\right)^{+\epsilon}
				\comCov{B_t}{X_t}_M
				\left(P^M_t\right)^{+\epsilon}\right]
				\mathrm{d}t
				\\&\phantom{=}+
				\mathrm{tr}\left[\left(P^M_t\right)^{+\epsilon}
				\left(P^M_t H^\mathrm{T}_t+\tilde{C}_t\Gamma^{\mathrm{T}}_t\right)
				R^{-1}_t
				\left(H_t+\Gamma_t\tilde{C}_t^{\mathrm{T}}\left(P^M_t\right)^{+}\right)
				P^M_t
				\left(P^M_t\right)^{+\epsilon}\right]
				\mathrm{d}t
				\\&\phantom{=}-
				\mathrm{tr}\left[\left(P^M_t\right)^{+\epsilon}
				\tilde{C}_t\tilde{C}^{\mathrm{T}}_t
				\left(P^M_t\right)^{+\epsilon}\right]
				\mathrm{d}t
				\\&\phantom{=}-
				\frac{1}{M} \sum_{i=1}^{M}
				\mathrm{tr}\left[\left(P^M_t\right)^{+\epsilon}
				C_t(X^i_t) C_t(X^i_t)^{\mathrm{T}}
				\left(P^M_t\right)^{+\epsilon}\right]
				\mathrm{d}t
				\\&\phantom{=}+	
				\frac{4}{(M-1)^2}\sum_{j=1}^{M}
				\mathcal{S}\left(\left(P^M_t\right)^{+\epsilon},X^j_t-x^M_t,C_t(X^j)\right)
				\mathrm{d}t
				\\&\phantom{=}+
				\frac{4}{(M-1)^2}\sum_{j=1}^{M}
				\mathcal{S}\left(\left(P^M_t\right)^{+\epsilon},X^j_t-x^M_t,\tilde{C}_t\right)
				\mathrm{d}t
				+
				\mathrm{d}\mathcal{M}^{\epsilon}_t.
			\end{split}
		\end{align}
		
		We aim to estimate this equation linearly and then use a Grönwall argument to bound $\mathrm{tr}\left(P^M\right)^{+\epsilon}$. To this end let us denote by $e_i$ the 
		$i$-th canonical basis vector of $\R^{d_x}$ and by $P^M_t=Q^{\mathrm{T}}\mathrm{diag}\left(\lambda^1_t,\cdots,\lambda^{d_x}_t\right)Q$ the singular value decomposition of $P^M$. Then we first note that for arbitrary $k$
		\begin{align*}
			&\sum_{j=1}^{M}\mathrm{tr}\left[
			\left(P^M_t\right)^{+\epsilon}
			\left(X^j_t-x^M_t\right)C(X^j)_{k}^{\mathrm{T}}
			\left(P^M_t\right)^{+\epsilon}
			\left(X^j_t-x^M_t\right)C(X^j)_{k}^{\mathrm{T}}
			\left(P^M_t\right)^{+\epsilon}
			\right] 
			\\&=
			\sum_{i=1}^{d_x} 
			\frac{
				\sum_{j=1}^{M} 
				\left(e^{\mathrm{T}}_i Q \left(X^j_t-x^M_t\right)\right)
				\left(C(X^j)_{k}^{\mathrm{T}}\left(P^M_t\right)^{+\epsilon}		\left(X^j_t-x^M_t\right)\right)
				\left(C(X^j)_{k}^{\mathrm{T}}Q^{\mathrm{T}} e_i\right)
			}{\left(\lambda^i_t+\epsilon\right)^2}
			\\&\leq
			\norm{C_k}{\infty}^2
			\sum_{i=1}^{d_x} \frac{1}{\left(\lambda^i_t+\epsilon\right)^2}
			\sum_{j=1}^{M} 
			\left|e^{\mathrm{T}}_i Q \left(X^j_t-x^M_t\right)\right|
			\left|\left(P^M_t\right)^{+\epsilon}\left(X^j_t-x^M_t\right)\right|
			\\&\leq
			\left(M-1\right)\norm{C_k}{\infty}^2
			\sum_{i=1}^{d_x} \frac{\sqrt{\lambda^i_t}}{\left(\lambda^i_t+\epsilon\right)^2}
			\sqrt{\sum_{l=1}^{d_x} \frac{\lambda^l_t}{\left(\lambda^l_t+\epsilon\right)^2}}.
		\end{align*}
		
		For the last inequality we used that
		\begin{align*}
			\frac{1}{M-1} \sum_{j=1}^{M} 
			\left|e^{\mathrm{T}}_i Q \left(X^j_t-x^M_t\right)\right|^2
			&=
			\lambda^i_t
			~\text{and that}~
			\\
			\frac{1}{M-1} \sum_{j=1}^{M} 
			\left|\left(P^M_t\right)^{+\epsilon}\left(X^j_t-x^M_t\right)\right|^2
			&=
			\sum_{l=1}^{d_x} \frac{\lambda^l_t}{\left(\lambda^l_t+\epsilon\right)^2}.	
		\end{align*}
		
		Using the Cauchy--Schwarz inequality two times one then can verify immediately that
		\begin{align*}
			\sum_{i=1}^{d_x} \frac{\sqrt{\lambda^i_t}}{\left(\lambda^i_t+\epsilon\right)^2}
			\sqrt{\sum_{l=1}^{d_x} \frac{\lambda^l_t}{\left(\lambda^l_t+\epsilon\right)^2}}
			\leq
			\sqrt{d_x}~
			\sum_{i=1}^{d_x} \frac{1}{\left(\lambda^i_t+\epsilon\right)^2}.
		\end{align*}
		
		With this and the estimate
		{\small\begin{align*}
				&\frac{\sum_{j=1}^{M}\mathrm{tr}\left[
					\left(P^M_t\right)^{+\epsilon}
					\left(X^j_t-x^M_t\right)C(X^j)_{k}^{\mathrm{T}}
					\left(P^M_t\right)^{+\epsilon}
					C(X^j)_{k}
					\left(X^j_t-x^M_t\right)^{\mathrm{T}}
					\left(P^M_t\right)^{+\epsilon}
					\right]}{M-1}
				\\&=
				\sum_{i=1}^{d_x} \frac{\sum_{j=1}^{M} 
					\left(e^{\mathrm{T}}_i Q \left(X^j_t-x^M_t\right)\right)^2
					\left(C(X^j)_{k}^{\mathrm{T}}\left(P^M_t\right)^{+\epsilon}	C(X^j)_{k}\right)
				}{\left(\lambda^i_t+\epsilon\right)^2(M-1)}
				\\&\leq
				\norm{C_k}{\infty}^2
				\underbrace{\left|\left(P^M_t\right)^{+\epsilon}\right|}_{\sum_{i=1}^{d_x} \frac{1}{\lambda^i_t+\epsilon}}
				\underbrace{\sum_{i=1}^{d_x} \frac{\lambda^i_t}{\left(\lambda^i_t+\epsilon\right)^2}}
				_{\sum_{i=1}^{d_x} \frac{1}{\lambda^i_t+\epsilon}}
				\leq
				\norm{C_k}{\infty}^2 \sum_{i=1}^{d_x}\frac{1}{\left(\lambda^i_t+\epsilon\right)^2}
		\end{align*}}%
		we can bound the Itô correction term
		{\small\begin{align}\label{bound ito correction emp cov}
				\begin{split}
					&\frac{4}{(M-1)^2}\sum_{j=1}^{M}
					\left(\mathcal{S}\left(\left(P^M_t\right)^{+\epsilon},X^j_t-x^M_t,C_t(X^j)\right)
					+
					\mathcal{S}\left(\left(P^M_t\right)^{+\epsilon},X^j_t-x^M_t,\tilde{C}_t\right)\right)
					\\&\leq
					\frac{2}{M-1}
					\left(1+\sqrt{d_x}\right) \left(\norm{C}{\infty}^2+\left|\tilde{C}\right|^2\right)
					\sum_{i=1}^{d_x}\frac{1}{\left(\lambda^i_t+\epsilon\right)^2}.
				\end{split}
		\end{align}}%
		
		Similarly one shows that
		\begin{align}
			\begin{split}
				&\mathrm{tr}\left[\left(P^M_t\right)^{+\epsilon}
				\left(P^M_t H^\mathrm{T}_t+\tilde{C}_t\Gamma^{\mathrm{T}}_t\right)
				R^{-1}_t
				\left(H_t+\Gamma_t\tilde{C}_t^{\mathrm{T}}\left(P^M_t\right)^{+}\right)
				P^M_t
				\left(P^M_t\right)^{+\epsilon}\right]
				\\&\leq
				d_x \lambda_{\max}\left(H^{\mathrm{T}}_t R^{-1}_t H_t\right)
				+2\left|H^{\mathrm{T}}_t R^{-1}_t \Gamma_t  \tilde{C}^{\mathrm{T}}_t\right|
				\mathrm{tr}\left[\left(P^M_t\right)^{+\epsilon}\right]
				\\&\phantom{=}+
				\sum_{i=1}^{d_x} \frac{1}{\left(\lambda^j_t+\epsilon\right)^2}
				e_i^{\mathrm{T}} Q \tilde{C}_t R^{-1}_t\tilde{C}_t^{\mathrm{T}} Q^{\mathrm{T}} e_i
			\end{split}
		\end{align}
		and
		\begin{align}
			\begin{split}
				&\mathrm{tr}\left[\left(P^M_t\right)^{+\epsilon}
				\tilde{C}\tilde{C}^{\mathrm{T}}
				\left(P^M_t\right)^{+\epsilon}\right]
				\\&\phantom{=}+
				\frac{M-1}{M} \frac{1}{M} \sum_{i=1}^{M}
				\mathrm{tr}\left[\left(P^M_t\right)^{+\epsilon}
				C(X^i_t) C(X^i_t)^{\mathrm{T}}
				\left(P^M_t\right)^{+\epsilon}\right]
				\\&\phantom{=}+
				\frac{1}{M^2 (M-1)}\sum_{i=1}^{M}\sum_{j\neq i}
				\mathrm{tr}\left[\left(P^M_t\right)^{+\epsilon}
				C(X^j_t) C(X^j_t)^{\mathrm{T}}
				\left(P^M_t\right)^{+\epsilon}\right]
				\\&\geq
				\sum_{i=1}^{d_x}
				\frac{e_i^{\mathrm{T}} Q_t \tilde{C}\tilde{C}^{\mathrm{T}} Q^{\mathrm{T}}_t e_i}{\left(\lambda^j_t+\epsilon\right)^2}
				+
				\inf_{x}\lambda_{\min}\left(C(x) C(x)^{\mathrm{T}}\right)~
				\sum_{j=1}^{d_x} \frac{1}{\left(\lambda^j_t+\epsilon\right)^2}.
			\end{split}
		\end{align}
		
		Combining these inequalities we derive the stochastic differential inequality
		\begin{align*}
			\mathrm{d}\mathrm{tr}\left(P^M_t\right)^{+\epsilon}
			&\leq
			-\mathrm{tr}\left[\left(P^M_t\right)^{+\epsilon}
			\comCov{B_t}{X_t}_M
			\left(P^M_t\right)^{+\epsilon}\right]
			\mathrm{d}t
			-
			\underline{\gamma}^M~
			\sum_{j=1}^{d_x} \frac{1}{\left(\lambda^j_t+\epsilon\right)^2}
			\mathrm{d}t
			\\&\phantom{=}+
			d_x \lambda_{\max}\left(H^{\mathrm{T}}_t R^{-1}_t H_t\right)\mathrm{d}t
			+2\left|H^{\mathrm{T}}_t R^{-1}_t \Gamma_t  \tilde{C}^{\mathrm{T}}_t\right|
			\mathrm{tr}\left[\left(P^M_t\right)^{+\epsilon}\right]\mathrm{d}t+
			\mathrm{d}\mathcal{M}^{\epsilon}_{t}.
		\end{align*}
		
		Using Cauchy--Schwarz and the Taylor inequality we derive for any arbitrary $\delta>0$
		\begin{align*}
			&-\mathrm{tr}\left[\left(P^M_t\right)^{+\epsilon}
			\comCov{B_t}{X_t}_M
			\left(P^M_t\right)^{+\epsilon}\right]
			\\&\leq
			\sum_{j=1}^{d_x}
			\frac{
				\sqrt{\frac{2}{M-1}\sum_{i=1}^{M}
					\left(e^{\mathrm{T}}_j Q\left(B_t\left(X^i_t\right)-b^M_t\right)\right)^2}
				\sqrt{\frac{2}{M-1}\sum_{i=1}^{M}
					\left(e^{\mathrm{T}}_jQ\left(X^i_t-x^M_t\right)  \right)^2}
			}{\left(\lambda^j_t+\epsilon\right)^2}   
			\\&\leq
			\sum_{j=1}^{d_x}
			\frac{2\sqrt{\lambda^j_t}\mathrm{Lip}(B)\sqrt{\mathrm{tr}P^M_t}}
			{\left(\lambda^j_t+\epsilon\right)^2}
			\leq
			\sum_{j=1}^{d_x}
			\frac{\frac{\mathrm{Lip(B)}^2 \mathrm{tr} P^M_t}{\delta} +\delta\lambda^j_t}
			{\left(\lambda^j_t+\epsilon\right)^2}
			\\&\leq
			\sum_{j=1}^{d_x}
			\frac{\mathrm{Lip(B)}^2 \mathrm{tr} P^M_t}
			{\delta\left(\lambda^j_t+\epsilon\right)^2}
			+
			\delta\mathrm{tr}\left(P^M_t\right)^{+\epsilon}.
		\end{align*}
		
		By setting $\delta:=\frac{\mathrm{Lip(B)}^2 \mathrm{tr}P^M_t}{\underline{\gamma}^M}$ we can thus bound the evolution of $\mathrm{tr}\left(P^M_t\right)^{+\epsilon}$ even further by
		\begin{align}\label{evolution bound inverse emp cov}
			\begin{split}
				\mathrm{d}\mathrm{tr}\left(P^M_t\right)^{+\epsilon}
				&\leq
				\left(\delta+2\left|H^{\mathrm{T}}_t R^{-1}_t  \Gamma_t \tilde{C}^{\mathrm{T}}_t\right|\right)\mathrm{tr}\left(P^M_t\right)^{+\epsilon}
				\mathrm{d}t
				\\&\phantom{=}+
				d_x \lambda_{\max}\left(H^{\mathrm{T}}_t R^{-1}_t H_t\right)\mathrm{d}t
				+
				\mathrm{d}\mathcal{M}^{\epsilon}_t.
			\end{split}
		\end{align}
		
		We set $\bar{\xi}^{\kappa}:=\inf\left\{~t\geq0~:~\mathrm{tr}P^M_t>\kappa~\right\}$. Then using the stochastic Grönwall lemma as found in \cite[Theorem 4]{Scheutzow} on the differential inequality \eqref{evolution bound inverse emp cov}  implies 
		\begin{align*}
			\mathbb{E}
			\left[\sup_{t\leq T\wedge \bar{\xi}^\kappa}\left(\mathrm{tr}\left(P^M_t\right)^{+\epsilon}\right)^{1/2}\right]
			&\leq
			\left(\pi+1\right)
			\exp\left(\int_{0}^{T} \left(\left|H^{\mathrm{T}}_s R^{-1}_s \Gamma_t \tilde{C}^{\mathrm{T}}_s\right|+\frac{\mathrm{Lip(B)}^2 \kappa}{2\gamma}\right)\mathrm{d}s\right)
			\\&\phantom{=} \times 
			\mathbb{E}\left[\left(\mathrm{tr}\left(P^M_0\right)^{+\epsilon}
			+
			d_x
			\int_{0}^{T}
			\lambda_{\max}\left(H^{\mathrm{T}}_s R^{-1}_s H_s\right)\mathrm{d}s
			\right)^{1/2}\right].
		\end{align*}
		We note that in particular this bound is independent of the local martingale $M^{\epsilon}$ and its quadratic variation,which is the key advantage of the stochastic Grönwall Lemma compared to standard bounds based on the Burkholder--Davis--Gundy inequality.\newline
		
		Letting $\epsilon$ go to zero thus shows that $\sup_{t\leq T\wedge \bar{\xi}^\kappa}\mathrm{tr}\left(P^M_t\right)^{+}$ can not blow up for arbitrary $\kappa>0$. Thus on the interval $[0,\bar{\xi})$ the empirical covariance matrix $P^M$ will always be regular, which also implies $\bar{\xi}\leq\underline{\xi}$. Therefore it is left to show that blow ups can not occur in order to prove the well posedness.\newline
		
		\textbf{Step 3:} The rest of the proof is similar to the one found in \cite{StannatLange} for the uncorrelated case.\newline
		
		We note that since $P^M_t H^\mathrm{T}_t R^{-1}_t H_t P^M_t$ is symmetric positive semidefinite and since $\left|\left(P^M_t\right)^{+} P^M_t\right|\leq \sqrt{d_x}$, we have
		\begin{align*}
			-
			&\mathrm{tr}\left[\left(P^M_t H^\mathrm{T}_t+\tilde{C}_t\Gamma^{\mathrm{T}}_t\right)
			R^{-1}_t \left(H_t+\Gamma_t\tilde{C}^{\mathrm{T}}_t\left(P^M_t\right)^{+}\right)
			P^M_t\right]
			\\&=
			-\mathrm{tr}\left[P^M_t H^\mathrm{T}_t R^{-1}_t H_t P^M_t\right]
			-\mathrm{tr}\left[P^M_t H^\mathrm{T}_t R^{-1}_t \Gamma_t \tilde{C}^{\mathrm{T}}_t\left(P^M_t\right)^{+} P^M_t\right]
			\\&\phantom{=}-
			\mathrm{tr}\left[\tilde{C}_t \Gamma^{\mathrm{T}}_t R^{-1}_t H_t P^M_t\right]
			-\mathrm{tr}\left[\tilde{C}_t \Gamma^{\mathrm{T}}_t R^{-1}_t \Gamma_t \tilde{C}^{\mathrm{T}}_t\left(P^M_t\right)^{+} P^M_t\right]
			\\&\leq
			2d_x \left|\tilde{C}_t \Gamma^{\mathrm{T}}_t R^{-1}_t H_t\right|~\mathrm{tr}\left[P^M_t\right]
			+
			d_x \left|\tilde{C}_t \Gamma^{\mathrm{T}}_t R^{-1}_t \Gamma_t \tilde{C}^{\mathrm{T}}_t\right|.
		\end{align*}
		
		This gives us the stochastic differential inequality 
		\begin{align}\label{inequality evolution of tr PM}
			\begin{split}
				\mathrm{d} \mathrm{tr} P^{M}_t
				&\leq
				2 \left(\mathrm{Lip}\left(B\right)+ d_x \left|\tilde{C}_t R^{-1}_t H_t\right|\right) \mathrm{tr}P^M_t\mathrm{d}t
				+
				2d_x \left|\tilde{C}_t \Gamma^{\mathrm{T}}_t R^{-1}_t \Gamma_t \tilde{C}^{\mathrm{T}}_t\right| \mathrm{d}t
				\\&\phantom{=}+
				\sup_{s\leq t} 
				\mathrm{tr} \tilde{C}_t \tilde{C}_t^{\mathrm{T}} \mathrm{d}t
				+
				\sup_{s\leq t}\norm{C_s}{\infty}^2 \mathrm{d}t
				+\mathrm{d}\mathfrak{lm}_t
			\end{split}
		\end{align}
		
		Using the stochastic Grönwall Lemma \cite[Theorem 4]{Scheutzow} on the differential inequality \eqref{inequality evolution of tr PM} we can thus bound the covariance matrix by
		\begin{align*}
			\mathbb{E}\left[
			\sup_{t\leq T\wedge\bar{\xi}}
			\sqrt{\mathrm{tr} P^{M}_t}
			\right]
			&\leq
			(\pi+1)
			\exp\left(\int_{0}^{T} 
			\mathrm{Lip}\left(B\right)+ d_x \left|\tilde{C}_s \Gamma^{\mathrm{T}}_s R^{-1}_s H_s\right|
			~\mathrm{d}s \right)
			\\&\phantom{=}
			\sqrt{	d_x T~\sup_{t\leq T} \left|\tilde{C}_t \Gamma^{\mathrm{T}}_t R^{-1}_t \Gamma_t \tilde{C}^{\mathrm{T}}_t\right|
				+
				\frac{T M}{M-1}
				\sup_{t\leq T} 
				\left(\left|\tilde{C}_t\right|^2+\norm{C_s}{\infty}^2\right)
			}.
		\end{align*}
		
		Due to the assumptions, the right-hand side of this inequality is finite. It is easy to see that on the set of paths where $P^M$ is bounded, the empirical mean evolves at most linearly and thus one can easily verify the bound
		\begin{align*}
			\mathbb{E}\left[\sup_{t\leq T\wedge\bar{\xi}} \sqrt{\left|x^M_t\right|}\right]<+\infty,
		\end{align*}
		which shows that $\mathbb{P}\left(\bar{\xi}<T\right)=0$. Thus explosion can not occur in the time interval $[0,T]$, which let's us use standard localization arguments (see for example \cite[Theorem 3.4 and its proof]{Mao}) to conclude the existence of a unique global solution to the EnKBF.\newline
	\end{proof}

	\begin{Rmk}
		{ We could drop the assumptions $\underline{\gamma}^M>0$ and $P^M_0$ being invertible in the uncorrelated case $\tilde{C}=0$, where well posedness  was proven already in \cite{StannatLange}.  We could also drop these assumptions if we replaced the Moore--Penrose inverse $\left(P^M_t\right)^{+}$ in \eqref{nonlinear Ensemble Kalman} by a regularized inverse of the form
			\begin{align*}
				\left(P^M_t\right)^{+\epsilon,n}:=\left(\left(P^M_t\right)^n+\epsilon I\right)^{-1} \left(P^M_t\right)^{n-1}.
			\end{align*}
			for arbitrary but fixed $\epsilon>0$ and $n\in\N$. For $n=1$ this is similar to Ensemble inflation, a well known technique in data assimilation \cite{LawKellyStuart}, whereas for $n=2$ this is a canonic approximation of the Pseudoinverse. Since this regularization depends Lipschitz-continuously on $P^M$ it would have then sufficed to only carry out steps 1 and 3 in the proof of theorem \ref{well posedness EnKBF theorem}.\\
			If one tried to prove the well posedness for \eqref{nonlinear Ensemble Kalman} without the regularization for small ensemble sizes $M\leq d_x$, i.e. without the assumption that $P^M_0$ is invertible, one would try to show that $P^M_t$ never changes its rank. This could be done similar to step 2 in the proof of theorem \ref{well posedness EnKBF theorem}. However one would need to show that $\left(P^M_t\right)^{+\epsilon,2}$, instead of just bounding $\left(P^M_t\right)^{+\epsilon,1}$ as it was done in the proof, as the latter one will always explode for non-invertible $P^M_t$ when $\epsilon\to 0$, whereas $\left(P^M_t\right)^{+\epsilon,2}$ only does this when the matrix rank changes.}
	\end{Rmk}
	
	Now we again shift our focus to the mean field limit of \eqref{nonlinear Ensemble Kalman}.

	\section{Propagation of Chaos}\label{section propagation of chaos}
	
	In this section we aim to prove that the particles defined by the regularized EnKBF \eqref{nonlinear Ensemble Kalman} indeed converge to the solution of \eqref{nonlinear mf ensemble kalman filter - true inverse} for $M\to\infty$. A proof for the uncorrelated case with time independent signal diffusion $C_t(x)=C,~\tilde{C}_t=0$ can be found in \cite{StannatLange}. The first results in the linear Gaussian setting were obtained in \cite{DelMoralTugaut}.\newline
	
	We consider $M$ independent copies $\bar{X}^i,~i=1,\cdots,M$ of the solution $\tilde{X}$ to \eqref{nonlinear mf ensemble kalman filter - true inverse} satisfying
	\begin{align}\label{independent mf copies}
		\begin{split}
			\mathrm{d}\bar{X}^i_t
			&=
			B_t(\bar{X}^i_t)\mathrm{d}t+
			C_t(\bar{X}^i_t)\mathrm{d}W^i_t+
			\tilde{C}_t\mathrm{d}V^i_t
			\\&\phantom{=}+
			\left(\bar{P}_t H^\mathrm{T}_t+\tilde{C}_t\Gamma^{\mathrm{T}}_t\right) R^{-1}_t
			\left(\mathrm{d}Y_t-\frac{H_t\left(\bar{X}^i_t+\bar{m}_t\right)}{2}\mathrm{d}t\right)
			\\&\phantom{=}-
			\frac{\bar{P}_t H^\mathrm{T}_t+\tilde{C}_t\Gamma^{\mathrm{T}}_t}{2} R^{-1}_t \Gamma_t\tilde{C}^{\mathrm{T}}_t \bar{P}^{+}_t\left(\bar{X}^i_t-\bar{m}_t\right) \mathrm{d}t.
		\end{split}
	\end{align}

	\begin{Thm}\label{thm prop of chaos}
		We make the same assumptions as in Theorem \ref{Existence and uniqueness for mf}. Define the error term $r^i_t:=X^i_t-\bar{X}^i_t$, then the following result holds
		\begin{align}\label{error particle approximation}
			\begin{split}
				\sup_{t\leq T} \frac{1}{M}\sum_{i=1}^{M} |r^i_t|^2\xrightarrow{M\to\infty} 0
			\end{split}
		\end{align}
		in probability (and almost surely with respect to $Y$).
	\end{Thm}
	\begin{proof}
		By Remark \ref{consistency of conditions} the condition $\gamma>0$ also implies $\gamma^M>\frac{\gamma}{2}>0$ for sufficiently large $M$. Thus we can assume in the following that $M$ is large enough so that the well posedness of the ensemble filter proven in Theorem \ref{well posedness EnKBF theorem} holds.\newline

		For the sake of brevity in formulas let us define 
		\begin{align}\label{short notation gainInflation term}
			\begin{split}
				\psi_t(P):=\frac{P H^{\mathrm{T}}_t+\tilde{C}_t\Gamma^{\mathrm{T}}_t}{2} R^{-1}_t \Gamma_t \tilde{C}^{\mathrm{T}}_t P^{+}.
			\end{split}
		\end{align}
		
		We note that since both  $X^i$ and $\bar{X}^i$ share the same initial condition. Then we have
		\begin{align*}
			r^i_t&=
			\int_{0}^{t}\left(B_s(X^i_s)-B_s(\bar{X}^i_s)\right)\mathrm{d}s
			+
			\int_{0}^{t}\left(C_s(X^i_s)-C_s(\bar{X}^i_s)\right)\mathrm{d}W^i_s
			\\&\phantom{=}+
			\int_{0}^{t}\left(P^M_s-\bar{P}_s\right)H^\mathrm{T}_s R^{-1}_s
			\left(\mathrm{d}Y_s-\frac{H_s(\bar{X}^i_s+\bar{m}_s)}{2}\mathrm{d}s\right)
			\\&\phantom{=}-
			\frac{1}{2}
			\int_{0}^{t}\left(\bar{P}_s H^\mathrm{T}_s+\tilde{C}_s\Gamma^{\mathrm{T}}_s\right) R^{-1}_s H_s\left(\left(X^i_t-\bar{X}^i_t\right)+\left(x^M_s-\bar{m}_s\right)\right)\mathrm{d}s
			\\&\phantom{=}+
			\int_{0}^{t}\psi_s(\bar{P}_s)\left(\left(X^i_t-\bar{X}^i_t\right)-\left(x^M_s-\bar{m}_s\right)\right)\mathrm{d}s
			\\&\phantom{=}-
			\int_{0}^{t}\left(\psi_s(P^M_s)-\psi_s(\bar{P}_s)\right)\left(X^i_t-x^M_s\right)\mathrm{d}s.
		\end{align*}

		To take care of the integral involving $Y$ one could later on use a change of measure to turn $Y$ into a Brownian motion. We will not do this and instead use the explicit representation $\mathrm{d}Y_s=H_s\left(X^{\mathrm{ref}}_s\right)\mathrm{d}s+\Gamma dV_s$ to rewrite the equation above into
		\begin{align*}
			r^i_t
			&=\int_{0}^{t}~\mathrm{d}\left(X^i_s-\bar{X}^i_s\right)
			\\&=
			\int_{0}^{t}\left(B_s(X^i_s)-B_s(\bar{X}^i_s)-\frac{1}{2}\bar{P}_s H^{\mathrm{T}}_s R^{-1}_s H_s\left(X^i_s-\bar{X}^i_s\right)\right)\mathrm{d}s
			\\&\phantom{=}+
			\int_{0}^{t}\left(C_s(X^i_s)-C_s(\bar{X}^i_s)\right)\mathrm{d}W^i_s
			+
			\int_{0}^{t}\left(P^M_s-\bar{P}_s\right)H^\mathrm{T}_s \Gamma_s\mathrm{d}V_s
			\\&\phantom{=}+
			\int_{0}^{t}\left(P^M_s-\bar{P}_s\right)H^\mathrm{T}_s R^{-1}_s H_s
			\left(X^\mathrm{ref}_s-
			\frac{\bar{X}^i_s+\bar{m}_s}{2}\mathrm{d}s\right)\mathrm{d}s
			\\&\phantom{=}-
			\frac{1}{2}
			\int_{0}^{t}\left(\bar{P}_s H^\mathrm{T}_s+\tilde{C}_s\Gamma^{\mathrm{T}}_s\right)R^{-1}_s H_s\left(x^M_s-\bar{m}_s\right)\mathrm{d}s
			\\&\phantom{=}+
			\int_{0}^{t}\psi(\bar{P}_s)\left(\left(X^i_t-\bar{X}^i_t\right)-\left(x^M_s-\bar{m}_s\right)\right)\mathrm{d}s
			\\&\phantom{=}-
			\int_{0}^{t}\left(\psi(P^M_s)-\psi(\bar{P}_s)\right)\left(X^i_s-x^M_s\right)
			\mathrm{d}s.
		\end{align*}

		Using Itô's rule, the Lipschitz properties of the coefficients, as well as the fact that $2 a\cdot b\leq |a|^2+|b|^2$, we get
		\begin{align*}
			\frac{1}{M}\sum_{i=1}^{M}\left|r^i_t\right|^2
			&\leq
			\int_{0}^{t}
			\mathcal{L}_s\frac{1}{M}\sum_{i=1}^{M}\left|r^i_s\right|^2	
			\mathrm{d}s+
			\frac{2}{M}\sum_{i=1}^{M}\int_{0}^{t}r^i_s\cdot\left(C_s(X^i_s)-C_s(\bar{X}^i_s)\right)\mathrm{d}W^i_s
			\\&\phantom{=}+
			\frac{2}{M}\sum_{i=1}^{M}\int_{0}^{t}r^i_s\cdot\left(P^M_s-\bar{P}_s\right)H^\mathrm{T}_s\Gamma_s\mathrm{d}V_s
			\\&\phantom{=}+
			\int_{0}^{t}\left|P^M_s-\bar{P}_s\right|^2 \left|H^\mathrm{T}_s R^{-1}_s H_s\right|^2
			\left(\left|X^\mathrm{ref}_s\right|^2
			+
			\frac{1}{2M}\sum_{i=1}^{M}\left|\bar{X}^i_s\right|^2+\frac{\left|\bar{m}_s\right|^2}{2}\right)\mathrm{d}s
			\\&\phantom{=}-
			\int_{0}^{t}
			\underbrace{\left(\frac{1}{M}\sum_{i=1}^{M} r^i_s\right)}_{=x^M_s-\bar{x}^M_s}
			\cdot\mathcal{A}_s\left(x^M_s-\bar{m}_s\right)\mathrm{d}s
			+
			\int_{0}^{t}
			\left|H^\mathrm{T}_s R^{-1}_s H_s\right|\left|P^M_s-\bar{P}_s\right|^2\mathrm{d}s
			\\&\phantom{=}+
			\int_{0}^{t}\left|\psi_s(P^M_s)-\psi_s(\bar{P}_s)\right|^2
			\underbrace{\left(\frac{1}{M}\sum_{i=1}^{M}\left|X^i_s-x^M_s\right|^2\right)}_{=\mathrm{tr}P^M_s}
			\mathrm{d}s
			,
		\end{align*}
		where we denote
		\begin{align}\label{definition L and A in prop}
			\begin{split}
				\mathcal{L}_s
				&:=2~\mathrm{Lip}(B)
				+\mathrm{Lip}(C)^2
				+\left|\bar{P}_s H^{\mathrm{T}}_s R^{-1}_s H_s\right|
				+2\left|\psi_s\left(\bar{P}_s\right)\right|
				+2
				\\
				\mathcal{A}_s
				&:=
				\bar{P}_s H^\mathrm{T}_s R^{-1}_s H_s + \tilde{C}_s \Gamma^{\mathrm{T}}_s R^{-1}_s H_s +2\psi_s\left(\bar{P}_s\right)
			\end{split}
		\end{align}
		for the sake of brevity.\newline

		First we note that
		\begin{align*}
			&-
			\left(x^M_s-\bar{x}^M_s\right)
			\cdot \left(\bar{P}_s H^\mathrm{T}_s R^{-1}_s H_s + \tilde{C}_s \Gamma^{\mathrm{T}}_s R^{-1}_s H_s+2\psi_s\left(\bar{P}_s\right)\right)\left(x^M_s-\bar{m}_s\right)
			\\&\leq
			\frac{1}{M}\sum_{i=1}^{M}\left|r^i_s\right|^2+
			\left|\bar{P}_s H^\mathrm{T}_s R^{-1}_s H_s + \tilde{C}_s \Gamma^{\mathrm{T}}_s R^{-1}_s H_s+2\psi_s\left(\bar{P}_s\right)\right|^2
			\left|\bar{x}^M_s-\bar{m}_s\right|^2
		\end{align*}
		and therefore if we denote
		\begin{align}\label{semimartingale in prop}
			\begin{split}
				\mathfrak{lm}^M_t
				&:=
				\frac{2}{M}\sum_{i=1}^{M}\int_{0}^{t}r^i_s\cdot\left(C_s(X^i_s)-C_s(\bar{X}^i_s)\right)\mathrm{d}W^i_s
				\\&\phantom{=}+
				\frac{2}{M}\sum_{i=1}^{M}\int_{0}^{t}r^i_s\cdot\left(P^M_s-\bar{P}_s\right)H^\mathrm{T}_s R^{-1}_s \Gamma_s\mathrm{d}V_s,
			\end{split}
		\end{align}
		we derive
		\begin{align*}
			\frac{1}{M}\sum_{i=1}^{M}\left|r^i_t\right|^2
			&\leq
			\int_{0}^{t}
			\left(\mathcal{L}_s
			+1\right) \frac{1}{M}\sum_{i=1}^{M}\left|r^i_t\right|^2	
			\mathrm{d}s
			\\&\phantom{=}+
			\int_{0}^{t}\left|P^M_s-\bar{P}_s\right|^2 \left|H^\mathrm{T}_s R^{-1}_s H_s\right|^2
			\left(\left|X^\mathrm{ref}_s\right|^2
			+
			\frac{1}{2M}\sum_{i=1}^{M}\left|\bar{X}^i_s\right|^2+\frac{\left|\bar{m}_s\right|^2}{2}\right)\mathrm{d}s
			\\&\phantom{=}+
			\int_{0}^{t}
			\left|\bar{P}_s H^\mathrm{T}_s R^{-1}_s H_s + \tilde{C}_s \Gamma^{\mathrm{T}}_s R^{-1}_s H_s+\psi_s\left(\bar{P}_s\right)\right|^2
			\left|\bar{x}^M_s-\bar{m}_s\right|^2
			\mathrm{d}s
			\\&\phantom{=}+
			\int_{0}^{t}\left|\psi_s(P^M_s)-\psi_s(\bar{P}_s)\right|^2~
			\mathrm{tr}P^M_s~
			\mathrm{d}s
			+\mathfrak{lm}^M_t.
		\end{align*}
		
		Clearly $\psi_s$ is locally Lipschitz if both arguments are invertible. Let us denote its Lipschitz constant on
		\begin{align}\label{definition set in prop}
			\begin{split}
				\mathfrak{S}_{\kappa}=\left\{~A\in\R^{d_x\times d_y}~|~\left|A\right|\leq \kappa~\text{and}~\left|A^{-1}\right|\leq \kappa~\right\}
			\end{split}
		\end{align}
		by $\mathrm{Lip}_{\mathrm{loc}}\left(\psi,\kappa\right)$. Since we have a priori bounds for $\bar{P}$, both from above and below, we can choose $\kappa$ large enough so that $\bar{P}_t\in\mathfrak{S}_{\kappa}$ for all times $t\leq T$. On the event that both $P^M$ also stays in the set $\mathfrak{S}_{\kappa}$ this gives us the bound
		\begin{align*}
			\left|\psi(P^M_s)-\psi(\bar{P}_s)\right|
			\leq
			\mathrm{Lip}_{\mathrm{loc}}\left(\psi_s,\kappa\right) 
			\left|P_s^{M}-\bar{P}_s\right|.
		\end{align*}
		
		Now we note that
		\begin{align*}
			\frac{1}{M}\sum_{i=1}^{M}\left|\bar{X}^i_s\right|^2
			&\leq
			\frac{2}{M}\sum_{i=1}^{M}\left|\bar{X}^i_s-\bar{x}^{m}_s\right|^2+2\left|\bar{x}^{m}_s\right|^2
			\leq
			\frac{2(M-1)}{M}\mathrm{tr}\bar{P}^M_s+2\left|\bar{x}^{m}_s\right|^2,
		\end{align*}
		and that
		\begin{align*}
			\left|P_s^{M}-\bar{P}_s\right|
			&\leq
			\left(\sqrt{\frac{M}{M-1}}+\frac{M}{M-1}\right)
			\sqrt{\mathrm{tr}P^M_s+\mathrm{tr}\bar{P}^M_s}
			\left(\frac{1}{M}\sum_{i=1}^{M}\left|r^i_s\right|^2\right)^{1/2}
			\\&\phantom{=}+\left|\bar{P}^M_s-\bar{P}_s\right|,
		\end{align*}
		where $\bar{x}^M$ and $\bar{P}^M$ denote the empirical mean and covariance matrix of the i.i.d. copies $\bar{X}^{i},~i=1,\cdots,M$.\newline

		Therefore we derive that on the event that $P^M$ stays in $\mathfrak{S}_{\kappa}$ we have
		\begin{align}
			\begin{split}
				\frac{1}{M}\sum_{i=1}^{M}\left|r^i_t\right|^2
				&\leq
				\int_{0}^{t}
				\mathcal{L}^1_s
				\frac{1}{M}\sum_{i=1}^{M}\left|r^i_s\right|^2	
				\mathrm{d}s
				+
				\int_{0}^{t}
				\mathcal{L}^2_s
				\left|\bar{P}^M_s-\bar{P}_s\right|^2
				\mathrm{d}s
				\\&\phantom{=}+
				\int_{0}^{t}
				\mathcal{L}^3_s
				\left|\bar{x}^M_s-\bar{m}_s\right|^2
				\mathrm{d}s+\mathfrak{lm}^M_t,
			\end{split}
		\end{align}
		where
		\begin{align}\label{definition constants in prop}
			\begin{split}
				\mathcal{L}^1_s&:=
				\mathcal{L}_s+1
				+
				\left(1+\mathrm{Lip}_{\mathrm{loc}}\left(\psi_s,\kappa\right) \right)^2
				\left(\sqrt{\frac{M}{M-1}}+\frac{M}{M-1}\right)^2\left(\mathrm{tr}P^M_s+\mathrm{tr}\bar{P}^M_s \right) \mathcal{L}^2_s
				\\\mathcal{L}^2_s&:=
				2 \left(\left|H^\mathrm{T}_s R^{-1}_s H_s\right|^2+1
				\right)
				\left(\left|X^\mathrm{ref}_s\right|^2+
				\frac{(M-1)}{M}\left|\bar{P}^M_s\right|^2+\left|\bar{x}^{m}_s\right|^2+
				\frac{\left|\bar{m}_s\right|^2}{2}+1\right)
				\\\mathcal{L}^3_s&:=
				\left|\bar{P}_s\right|^2\left|H^\mathrm{T}_s R^{-1}_s H_s\right|^2.
			\end{split}
		\end{align}
		
		For arbitrary constant $\kappa>0$ we define the stopping time
		\begin{align}\label{definition stopping time in prop}
			\begin{split}
				\zeta_{\kappa}:=\inf\left\{
				~
				t\geq 0
				~:~
				P^M_t\not\in\mathfrak{S}_{\kappa}
				~\text{or}~
				\mathrm{tr}\bar{P}^{M}_t\leq\kappa
				~
				\right\}
			\end{split}
		\end{align}
		Then on the stochastic interval $[0,\zeta_{\kappa}]$ we can bound the three constants $\mathcal{L}^1_t,\mathcal{L}^2_t, \mathcal{L}^3_t$ uniformly by some $\mathcal{C}_{\kappa}$. Thus the stochastic Grönwall Lemma implies that there exists a constant $c>0$ such that
		\begin{align*}
			&\mathbb{E}\left[\sup_{t\leq T\wedge\zeta_{\kappa}}\sqrt{\frac{1}{M}\sum_{i=1}^{M}\left|r^i_t\right|^2}\right]
			\leq
			c
			\sqrt{\mathcal{C}_{\kappa}}\exp\left(\frac{T \mathcal{C}_{\kappa} }{2}\right)
			\mathbb{E}\left[
			\sqrt{\left|\bar{x}^M_s-\bar{m}_s\right|^2+\left|\bar{P}^M_s-\bar{P}_s\right|^2}
			\right]
		\end{align*}
		
		For fixed $\kappa$ the right hand side will converge to zero for $M\to+\infty$ due to the law of large numbers. Thus for every fixed $\kappa$ the error term converges in probability against zero. Since $\zeta_{\kappa}\to+\infty$ for $\kappa\to+\infty$. This concludes the proof.
		
	\end{proof}

	Finally we want to note that both the well posedness result in Theorem \ref{well posedness EnKBF theorem} and the propagation of chaos result in Theorem \ref{thm prop of chaos} can easily be adapted to the particle systems approximating the transport based mean field EnKBF \eqref{mf limit of transport EnKBF} and the vanilla EnKBF \eqref{vanilla EnKBF}.\newline

	\textbf{Funding}
	Sebastian Ertel is supported by
	Deutsche Forschungsgemeinschaft through 
	\emph{IRTG 2544 - Stochastic Analysis in Interaction}.\\	
	The research of Wilhelm Stannat has been partially funded by Deutsche Forschungsgemeinschaft (DFG) - Project-ID 318763901 - SFB 1294.

	\bibliographystyle{alpha} 

\end{document}